\documentclass[a4paper,german,final,twoside,notitlepage,12pt]{article}

\usepackage[T1]{fontenc}
\usepackage{amsfonts}
\usepackage{amssymb}
\usepackage{amsthm}
\usepackage{amsmath}
\usepackage[all]{xy}
\usepackage{color}
\usepackage{graphicx}
\usepackage{hyperref}

\newtheorem{definition}{Definition}
\newtheorem{corollary}{Corollary}
\newtheorem{proposition}{Proposition}
\newtheorem{lemma}{Lemma}
\newtheorem{theorem}{Theorem}

\def\C {\mathbb{C}}

\def\id {\mathrm{id}}
\def\lto{\longrightarrow}
\def\lmaps{\;:\;}

\makeatletter
\def\proof{{Proof. }}
\def\endofproof{\hfill$\square$\\}
\def\adress#1{\gdef\@adress{#1}}
\def\@adress{}
\def\preprint#1{\gdef\@preprint{#1}}
\def\@preprint{}
\def\@maketitle{
  \newpage
  \noindent
  \begin{tabular}{cc}
    \begin{minipage}[c]{0.4\textwidth}
      \begin{flushleft}
        \includegraphics[width=110pt]{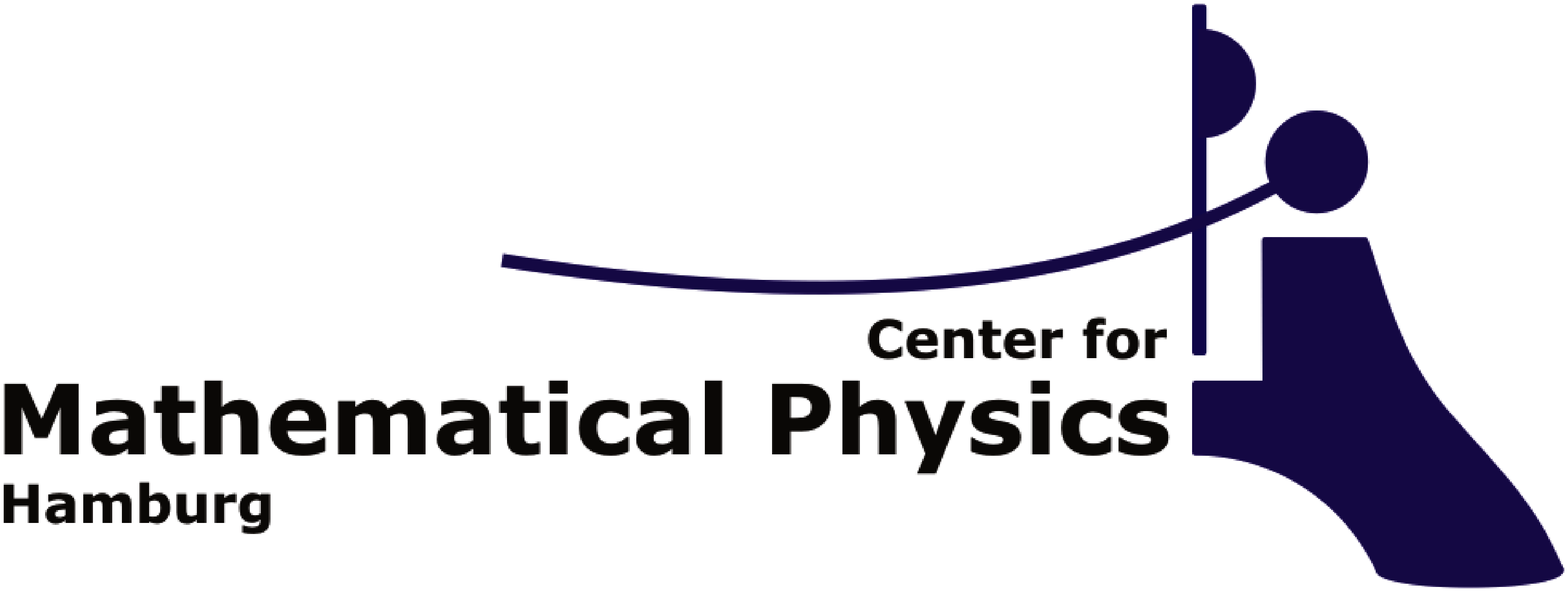}
      \end{flushleft}  
    \end{minipage}&
    \begin{minipage}[c]{0.55\textwidth}
      \begin{flushright}
      {\small\sf\@preprint}
      \end{flushright}
    \end{minipage}
  \end{tabular}
  \vskip 3cm
  \begin{center}
    \LARGE\@title
    \if!\@author!\else \vskip 0.5cm \large\@author\fi
    \if!\@adress!\else \vskip 0.5cm \normalsize\@adress\fi
  \end{center}
  \vskip 2cm
}
\newcommand{\alxy}[1]{\begin{aligned}\xymatrix{#1}\end{aligned}}
\newcommand{\alxydim}[2]{\begin{aligned}\xymatrix#1{#2}\end{aligned}}

\makeatother

\begin{document}

\title{More Morphisms between Bundle Gerbes}
\author{Konrad Waldorf}

\adress{Department Mathematik\\Schwerpunkt Algebra und Zahlentheorie\\Universit\"at
Hamburg\\Bundesstra\ss e 55\\D--20146 Hamburg}

\preprint{arXiv:math.CT/0702652\\
Hamburger Beiträge zur Mathematik Nr. 266\\
ZMP-HH/07-01}

\maketitle

\thispagestyle{empty}

\begin{abstract}
Usually bundle gerbes are considered as objects of a 2-groupoid, whose
 1-morphisms, called stable isomorphisms, are all invertible.  I introduce
 new 1-morphisms
 which
include stable isomorphisms, trivializations and bundle
gerbe modules. They fit into the structure of a 2-category of bundle
gerbes, and lead to natural  definitions of surface holonomy for closed surfaces, surfaces
with boundary, and unoriented closed surfaces. 
\end{abstract}

\newpage

\thispagestyle{empty}

\tableofcontents

\newpage

\setcounter{page}{1}

\section*{Introduction}
\addcontentsline{toc}{section}{Introduction}

From several perspectives it becomes clear that bundle gerbes are objects
in a 2-category: from the bird's-eye view of algebraic geometry, where gerbes appear as some kind of stack, or in topology, where they appear as one possible categorification
of a line bundle, but also from a worm's-eye view on the definitions of bundle
gerbes and their morphisms, which show that there have to be morphisms between the morphisms.

In \cite{stevenson1} a 2-groupoid is defined, whose objects are bundle gerbes,
and whose 1-morphisms are stable isomorphisms. To explain a few details,
recall that  bundle gerbes
are defined using surjective submersions $\pi:Y \to M$, and that a stable 
isomorphism $\mathcal{A}:\mathcal{G}_1 \to \mathcal{G}_2$ between two bundle
gerbes $\mathcal{G}_1$ and $\mathcal{G}_2$ with surjective submersions $\pi_1:Y_1
\to M$ and $\pi_2:Y \to M$ consists of a certain  line bundle $A$ over the fibre product
$Y_1 \times_M Y_2$. 2-morphisms between stable isomorphisms are  morphisms $\beta:A \to A'$ of those line bundles, obeying a compatibility constraint. Many examples of surjective submersions
arise from open covers $\lbrace U_{\alpha} \rbrace_{\alpha\in A}$ of $M$
by taking $Y$ to be the disjoint union of the open sets $U_{\alpha}$ and
$\pi$ to be the projection $(x,\alpha) \mapsto x$.  From this point of view, fibre products
$Y_1 \times_M Y_2$ correspond the common refinement of two open covers. So, the line bundle $A$ of a stable isomorphism
lives over the common refinement of the open covers of the two bundle gerbes.

Difficulties with this definition of stable isomorphisms arise when
two stable isomorphisms $\mathcal{A}:\mathcal{G}_1 \to \mathcal{G}_2$ and
$\mathcal{A}':\mathcal{G}_2 \to \mathcal{G}_3$ are going to be composed:
one has to define a line bundle $\tilde A$ over $Y_1 \times_M Y_3$ using the line bundles
$A$ over $Y_1 \times_M Y_2$ and $A'$ over $Y_2 \times_M Y_3$. In
\cite{stevenson1} this problem is solved using
 descent theory for line bundles.    
 
 In this note, I present another definition of 1-morphisms between bundle
gerbes (Definition \ref{def1}). Compared to stable isomorphisms, their definition is relaxed in two aspects: 
\begin{enumerate}
\item[1)] 
the line bundle is replaced by a certain vector bundle of rank possibly higher than 1. 
\item[2)]
this vector bundle is defined over a smooth manifold $Z$ with
surjective submersion $\zeta: Z \to Y_1 \times_M Y_2$. In terms of open covers,
the vector bundle lives over a refinement of the common refinement of the
open covers of the two bundle gerbes.
\end{enumerate}
Stable isomorphisms appear as a particular case of this relaxed definition.
I also give a generalized  definition
of 2-morphisms between such 1-morphisms (Definition \ref{def2}). 
Two goals are achieved by  this new type of morphisms between bundle
gerbes. First, relaxation 1) produces many 1-morphisms which are not invertible,
in contrast to the stable isomorphisms in \cite{stevenson1}. To be more precise, a 1-morphism is invertible
if and only if its vector bundle has rank 1 (Proposition \ref{prop2}). The
non-invertible  1-morphisms provide a new  formulation of left and right  bundle gerbe modules (Definition \ref{def4}). Second, relaxation 1) erases
 the difficulties with the composition of 1-morphisms: the vector bundle $\tilde A$ of the composition $\mathcal{A}'\circ \mathcal{A}$
 is just defined to be the tensor product
 of the vector bundles $A$ and $A'$ of the two 1-morphisms pulled back to $Z
 \times_{Y_2} Z'$, where $A$ over $Z$ is the vector bundle of $\mathcal{A}$
 and $A'$ over $Z'$ is the vector bundle of $\mathcal{A}'$. The composition defined like that is strictly associative (Proposition \ref{lem2}). This way we
end up with a strictly associative 2-category $\mathfrak{BGrb}(M)$ of bundle gerbes over
$M$. The aim of this note is to show that a good understanding of this 2-category
can be useful.    

This note is organized as follows. Section \ref{sec3} contains the 
definitions and properties  of the  2-category $\mathfrak{BGrb}(M)$ of bundle
gerbes over $M$. We also equip this 2-category with a monoidal structure,
pullbacks and a duality. Section \ref{sec5} relates our new definition of 1-morphisms between bundle
gerbes to the one of a stable isomorphism: two bundle gerbes
are isomorphic objects in $\mathfrak{BGrb}(M)$ if and only if they are stably
isomorphic (Corollary \ref{co1}). 
In section \ref{sec4} we present a unified view on important
structure related to bundle gerbes in terms of the new morphisms of the 2-category $\mathfrak{BGrb}(M)$: 
\begin{itemize}

\item[a)]
a \textit{trivialization} of a bundle gerbe $\mathcal{G}$ is a 1-isomorphism $\mathcal{A}:\mathcal{G}
\to \mathcal{I_{\rho}}$ from $\mathcal{G}$ to a trivial bundle gerbe
$\mathcal{I}_{\rho}$ given by a 2-form $\rho$ on $M$.

\item[b)]
a \textit{bundle gerbe module} of a bundle gerbe $\mathcal{G}$ is a (not
necessarily invertible) 1-morphism $\mathcal{E}:\mathcal{G}
\to \mathcal{I}_{\omega}$ from $\mathcal{G}$ to a trivial bundle gerbe $\mathcal{I}_{\omega}$.

\item[c)]
a \textit{Jandl structure} on a bundle gerbe $\mathcal{G}$ over $M$ is
a triple $(k,\mathcal{A},\varphi)$ of an involution $k$ of $M$,
a 1-isomorphism $\mathcal{A}:k^{*}\mathcal{G} \to \mathcal{G}^{*}$ and a certain 2-morphism $\varphi:k^{*}\mathcal{A} \Rightarrow \mathcal{A}^{*}$.

\end{itemize} 
Then we demonstrate how this understanding in combination with the
properties  of the
2-category $\mathfrak{BGrb}(M)$ can  be employed to give convenient  definitions of
surface holonomy. For this purpose we classify  the
morphisms between trivial bundle gerbes:
there is an equivalence of categories
\begin{equation*}
\mathfrak{Hom}(\mathcal{I}_{\rho_1},\mathcal{I}_{\rho_2}) \cong \mathfrak{Bun}_{\rho_2-\rho_1}(M)
\end{equation*}
between the morphism category between the trivial bundle gerbes  $\mathcal{I}_{\rho_1}$ and $\mathcal{I}_{\rho_2}$ and the category of vector bundles over $M$
for which the trace of the curvature gives the 2-form $\rho_2-\rho_1$ times
its rank. 

The interpretation of bundle gerbe modules and Jandl structures in terms
of
 morphisms between bundle gerbes is one step to understand the relation between two approaches to two-dimensional conformal field theories: on the one hand the Lagrangian
approach with a metric and a bundle gerbe $\mathcal{G}$ being the  relevant structure \cite{gawedzki1}
and
on the other
hand the algebraic approach in which a special symmetric Frobenius algebra
object $A$
in a modular tensor category $\mathcal{C}$ plays this role \cite{fuchs}.
Similarly as bundle gerbes, special symmetric Frobenius algebra objects in $\mathcal{C}$ form a 2-category, called $\mathcal{F}rob_{\mathcal{C}}$. In both approaches it is well-known how boundary conditions
have to be imposed. In the Lagrangian approach one chooses a D-brane: a submanifold
$Q$ of the target space together with a bundle gerbe module for the bundle gerbe $\mathcal{G}$ restricted to $Q$ \cite{gawedzki4}. In the algebraic approach one chooses
a 1-morphism from $A$ to the tensor unit $I$ of $\mathcal{C}$ (which is trivially
a special symmetric Frobenius algebra object) in the 2-category $\mathcal{F}rob_{\mathcal{C}}$
\cite{schweigert1}.
Now that we understand a gerbe module as a 1-morphism from $\mathcal{G}$ to $\mathcal{I}_{\omega}$ we have found a common principle in both approaches.
A similar success is made for unoriented conformal field theories. In the
Lagrangian approach, the bundle gerbe $\mathcal{G}$ has to be endowed with
a Jandl structure \cite{schreiber1}, which is in particular a 1-isomorphism
$k^{*}\mathcal{G} \to \mathcal{G}^{*}$ to the dual bundle gerbe $\mathcal{G}^{*}$.
In the algebraic approach one has to choose a certain algebra isomorphism $A\to A^{\mathrm{op}}$
from $A$ to the opposed algebra \cite{fuchs2}.

\paragraph{Acknowledgements.}
I would like to thank Christoph Schweigert for his advice and encouragement,
and Urs Schreiber for the many helpful discussions on 2-categories.

\paragraph{Conventions.} 
Let us fix the following conventions for the whole
article: by \textit{vector bundle} I refer to a complex vector bundle of
finite rank,
equipped with a hermitian structure and with a connection respecting this
hermitian structure. Accordingly, a \textit{morphism of vector bundles} is
supposed to respect both the hermitian structures and the connections. In particular, a \textit{line bundle} is a vector bundle in the above sense of rank one. The symmetric monoidal category $\mathfrak{Bun}(M)$, which is
formed by all  vector bundles over a smooth manifold $M$  and their morphisms
in the above sense,
is  for simplicity  tacitly replaced by an equivalent strict tensor category.

\section{The 2-Category of Bundle Gerbes}

\label{sec3}

Summarizing, the  2-category $\mathfrak{BGrb}(M)$ of bundle gerbes over a smooth manifold $M$ consists of
the following structure:\begin{enumerate}
\item 
A class of objects -- bundle gerbes over $M$. 
\item
A  morphism category $\mathfrak{Hom}(\mathcal{G},\mathcal{H})$   for each pair $\mathcal{G}$,
$\mathcal{H}$  of bundle gerbes, whose objects are called 1-morphisms and
are denoted by $\mathcal{A}:\mathcal{G} \to \mathcal{H}$, and whose morphisms
are called 2-morphisms and are denoted $\beta:\mathcal{A}_1 \Rightarrow \mathcal{A}_2$.
\item 
A  composition functor
\begin{equation*}
\circ \lmaps \mathfrak{Hom}(\mathcal{H},\mathcal{K}) \times \mathfrak{Hom}(\mathcal{G},\mathcal{H})
\lto \mathfrak{Hom}(\mathcal{G},\mathcal{K})
\end{equation*}
for each triple $\mathcal{G},\mathcal{H},\mathcal{K}$ of  bundle gerbes.
\item
An identity 1-morphism $\id_{\mathcal{G}}:\mathcal{G} \to \mathcal{G}$ for each bundle gerbe $\mathcal{G}$ together with natural 2-isomorphisms
\begin{equation*}
\rho_{\mathcal{A}}\lmaps \id_{\mathcal{H}} \circ \mathcal{A}\Longrightarrow\mathcal{A}
\quad\text{ and }\quad
\lambda_{\mathcal{A}}\lmaps \mathcal{A}\circ \id_{\mathcal{G}}\Longrightarrow \mathcal{A}
\end{equation*}
associated to every 1-morphism $\mathcal{A}:\mathcal{G} \to \mathcal{H}$.
\end{enumerate}
 This structure satisfies
the axioms of a strictly associative 2-category:
\begin{list}{}{\leftmargin=1.3cm\labelwidth=1.1cm\labelsep=0.2cm}
\item[\hypertarget{2C1}{(2C1)}] 
For  three 1-morphisms $\mathcal{A}:\mathcal{G}_1 \to \mathcal{G}_2$,
$\mathcal{A}':\mathcal{G}_2 \to \mathcal{G}_3$ and $\mathcal{A}'':\mathcal{G}_3
\to \mathcal{G}_4$, the composition functor satisfies
\begin{equation*}
\mathcal{A}'' \circ (\mathcal{A}' \circ \mathcal{A}) \;=\; (\mathcal{A}'' \circ \mathcal{A}') \circ \mathcal{A}\text{.}
\end{equation*}

\item[\hypertarget{2C2}{(2C2)}]
For 1-morphisms $\mathcal{A}:\mathcal{G}_1 \to \mathcal{G}_2$ and $\mathcal{A}':\mathcal{G}_2
\to \mathcal{G}_3$, the 2-isomorphisms $\lambda_{\mathcal{A}}$ and $\rho_{\mathcal{A}}$ satisfy the equality
\begin{equation*}
\id_{\mathcal{A}'} \circ \rho_{\mathcal{A}} \;=\; \lambda_{\mathcal{A}'} \circ
\id_{\mathcal{A}}
\end{equation*}
as 2-morphisms from $\mathcal{A}' \circ \id_{\mathcal{G}_2} \circ \mathcal{A}$
to $\mathcal{A}'\circ \mathcal{A}$.
\end{list}

\medskip

The following two subsections contain the definitions of the structure of
the 2-category $\mathfrak{BGrb}(M)$. 
The two  axioms are proved in Propositions  \ref{lem2} and \ref{lem1}. The
reader who is not interested in these details may directly continue   with section \ref{sec4}.

\subsection{Objects and Morphisms}

\label{sec1}

The definition of the objects of the 2-category $\mathfrak{BGrb}(M)$ -- the bundle gerbes over
$M$ -- is the usual one, just like, for instance, in \cite{murray,stevenson1,gawedzki1}.
Given a  surjective submersion $\pi:Y \to M$  we use the notation $Y^{[k]} := Y \times_M ... \times_M Y$ for the $k$-fold fibre product,
which is again a smooth manifold. Here we consider fibre products to be strictly
associative for simplicity. For the canonical projections between fibre products we use the notation $\pi_{i_1...i_k}:Y^{[n]} \to Y^{[k]}$. 

\begin{definition}
\label{def_gerbe}
A bundle gerbe $\mathcal{G}$ over
a smooth manifold
$M$ consists
of the following structure:
\begin{enumerate}

\item
a surjective submersion $\pi: Y
\to M$,

\item
a line bundle $L$ over $Y^{[2]}$,

\item
a 2-form $C \in \Omega^2(Y)$, and

\item
an isomorphism
\begin{equation*}
\mu \lmaps \pi_{12}^{*}L \otimes \pi_{23}^{*}L
\lto \pi_{13}^{*}L
\end{equation*}
of line bundles over $Y^{[3]}$.
\end{enumerate}
This structure has to satisfy two axioms:
\begin{list}{}{\leftmargin=1.3cm\labelwidth=1.1cm\labelsep=0.2cm}

\item[\hypertarget{G1}{(G1)}]
The curvature of $L$ is fixed by
\begin{equation*}
\mathrm{curv}(L) = \pi_2^{*}C -
\pi_1^{*}C\text{.}
\end{equation*}

\item[\hypertarget{G2}{(G2)}]
$\mu$ is associative
in the sense that the diagram
\begin{equation*}
\alxydim{@C=3cm@R=1.5cm}{\pi_{12}^{*}L \otimes \pi_{23}^{*}L \otimes \pi_{34}^{*}L \ar[r]^-{\pi_{123}^{*}\mu\otimes
\id} \ar[d]_{\id \otimes \pi_{234}^{*}\mu} & \pi_{13}^{*}L \otimes \pi_{34}^{*}L
\ar[d]^{\pi_{134}^{*}\mu} \\ \pi_{12}^{*}L \otimes \pi_{24}^{*}L \ar[r]_-{\pi_{124}^{*}
\mu} & \pi_{14}^{*}L}
\end{equation*}
of isomorphisms of  line bundles
  over $Y^{[4]}$ is commutative. 
\end{list}
\end{definition}

To give an example of a bundle gerbe, we
introduce trivial bundle gerbes. Just as
for every 1-form $A \in \Omega^1(M)$ there
is the (topologically) trivial line bundle over $M$ having
this 1-form as its connection 1-form, we find a
trivial bundle gerbe for every 2-form $\rho\in\Omega^2(M)$.
Its surjective submersion is  the identity $\id: M \to M$, and its 2-form
is $\rho$. Its line bundle over $M \times_M M\cong M$ is the trivial line bundle with the trivial
connection, and its isomorphism
is the identity isomorphism between trivial
line bundles. Now, axiom \hyperlink{G1}{(G1)} 
is satisfied
since $\mathrm{curv}(L)=0$ and
$\pi_1 = \pi_2 = \id_{M}$. The
associativity axiom \hyperlink{G2}{(G2)} is surely
satisfied by the identity isomorphism.
Thus we have defined a bundle gerbe,
which we denote by $\mathcal{I}_{\rho}$. 

\medskip

It should not be unmentioned that the geometric nature of bundle gerbes allows
explicit constructions of all (bi-invariant) bundle gerbes over all compact, connected
and simple Lie groups \cite{gawedzki1,meinrenken1,gawedzki2}. It becomes
in particular essential that a surjective submersion $\pi:Y \to M$ is more
general than an open cover of $M$.

\medskip

An important consequence of the existence of the isomorphism $\mu$ in the
structure of a bundle gerbe $\mathcal{G}$ is that the line bundle $L$ restricted
to the image of the diagonal embedding $\Delta:Y \to Y^{[2]}$ is canonically
trivializable (as a line bundle with connection):

\begin{lemma}
\label{lem4}
There is a canonical isomorphism $t_{\mu}: \Delta^{*}L \to 1$ of line bundles
over $Y$, which satisfies
\begin{equation*}
\pi_1^{*}t_{\mu} \otimes \id = \Delta_{112}^{*}\mu 
\quad\text{ and }\quad
\id\otimes \pi_2^{*}t_{\mu} = \Delta_{122}^{*}\mu 
\end{equation*}
as isomorphisms of line bundles over $Y^{[2]}$, where $\Delta_{112}:Y^{[2]} \to Y^{[3]}$ duplicates the first and $\Delta_{122}:Y^{[2]} \to Y^{[3]}$
duplicates the second factor. 
\end{lemma}

\proof
The isomorphism $t_{\mu}$ is defined using the canonical pairing with the
dual line bundle $L^{*}$ (which is strict by convention) and the isomorphism $\mu$:
\begin{equation}
\alxydim{@C=2cm}{\Delta^{*}L = \Delta^{*}L \otimes \Delta^{*}L \otimes \Delta^{*}L^{*} \ar[r]^-{\Delta^{*}\mu \otimes \id} &
\Delta^{*}L \otimes \Delta^{*}L^{*}=1}
\end{equation}
The two claimed equations follow from the associativity axiom \hyperlink{G2}{(G2)} by pullback
of the diagram along $\Delta_{1222}$ and $\Delta_{1112}$ respectively. 
\endofproof

\medskip

Now we  define the  category $\mathfrak{Hom}(\mathcal{G}_1,\mathcal{G}_2)$
for two bundle gerbes $\mathcal{G}_1$ and $\mathcal{G}_2$, to whose
 structure we refer by the same letters as in Definition \ref{def_gerbe} but with indices 1 or
2 respectively. 

\begin{definition}
\label{def1}A 1-morphism $\mathcal{A}:\mathcal{G}_1
\to \mathcal{G}_2$  consists of the following structure:
\begin{enumerate}
\item
a surjective submersion $\zeta:Z \to Y_1 \times_M Y_2$,

\item
a vector bundle $A$ over $Z$, and

\item
an isomorphism
\begin{equation}
\label{2}
\alpha \lmaps L_1  \otimes
\zeta_{2}^* A \lto  \zeta_{1}^* A \otimes L_2
\end{equation}
of vector bundles over $Z \times_M Z$.

\end{enumerate}
This structure has to satisfy two axioms:
\begin{list}{}{\leftmargin=1.5cm\labelwidth=1.3cm\labelsep=0.2cm}
\item[\hypertarget{1M1}{(1M1)}] 
The curvature of $A$ obeys
\begin{equation*}
\frac{1}{n}\mathrm{tr}(\mathrm{curv}(A))=C_2 - C_1\text{,}
\end{equation*}
where $n$ is the rank of the vector bundle $A$.

\item[\hypertarget{1M2}{(1M2)}]
 The isomorphism $\alpha$ 
is compatible with the isomorphisms $\mu_1$ and
$\mu_2$ of the gerbes $\mathcal{G}_1$ and $\mathcal{G}_2$ in the sense
that the diagram
\begin{equation*}
\alxydim{@C=2.5cm}{\zeta_{12}^{*}L_1 \otimes \zeta_{23}^{*}L_1 \otimes \zeta_3^{*}A \ar[r]^-{\mu_1
\otimes \id} \ar[d]_{\id \otimes \zeta_{23}^{*}\alpha} & \zeta_{13}^{*}L_1
\otimes \zeta_3^{*}A \ar[dd]^{\zeta_{13}^{*}\alpha} \\ \zeta_{12}^{*}L_1
\otimes \zeta_2^{*}A \otimes \zeta_{23}^{*}L_2 \ar[d]_{\zeta_{12}^{*}\alpha
\otimes \id} & \\ \zeta_1^{*}A \otimes \zeta_{12}^{*}L_2 \otimes \zeta_{23}^{*}L_2
\ar[r]_-{\id \otimes \mu_2} & \zeta_1^{*}A \otimes \zeta_{13}^{*}L_2}
\label{1}
\end{equation*}
of isomorphisms of vector bundles
over $Z\times_M Z \times_M Z$ is commutative.
\end{list}
\end{definition}

Here we work with the following simplifying notation: we have not introduced
notation for
the canonical projections 
$Z \to Y_1$ and $Z \to Y_2$, accordingly we don't write pullbacks with these maps. So in (\ref{2}), where the line bundles $L_i$ are pulled back along the induced map
$Z^{[2]} \to Y_i^{[2]}$ for $i=1,2$ and also in axiom \hyperlink{1M1}{(1M1)} which
is an equation of 2-forms on $Z$.

\medskip

As an example of a 1-morphism, we define the
 identity 1-morphism
\begin{equation} 
 \id_{\mathcal{G}}\lmaps\mathcal{G}
\lto \mathcal{G}
\end{equation}
 of a bundle gerbe $\mathcal{G}$ over $M$. 
It is defined by $Z:=Y^{[2]}$, the identity $\zeta:=\id_Z$, the line bundle $L$ of $\mathcal{G}$ over $Z$ and the isomorphism $\lambda$
defined by
\begin{equation}
\alxydim{@C=2cm}{\pi_{13}^{*}L \otimes \pi_{34}^{*}L \ar[r]^-{\pi_{134}^{*}\mu} & \pi_{14}^{*}L
\ar[r]^-{\pi_{124}^{*}\mu^{-1}} & \pi_{12}^{*}L \otimes \pi_{24}^{*}L\text{,}}
\end{equation}
where we identified $Z^{[2]} = Y^{[4]}$, $\zeta_2 = \pi_{34}$ and $\zeta_1=\pi_{12}$.
Axiom \hyperlink{1M1}{(1M1)} is the same as axiom \hyperlink{G1}{(G1)} for
the bundle gerbe $\mathcal{G}$ and axiom \hyperlink{1M2}{(1M2)} follows from  axiom \hyperlink{G2}{(G2)}.
 
\medskip

The following lemma introduces an important isomorphism of vector bundles
associated to every 1-morphism,
which will be used frequently in the definition of the structure of $\mathfrak{BGrb}(M)$
and also in section \ref{sec5}.  
\begin{lemma}
\label{lem3}
For any 1-morphism $\mathcal{A}:\mathcal{G}_1 \to \mathcal{G}_2$
there is a canonical isomorphism
\begin{equation*}
\mathrm{d}_\mathcal{A}\lmaps \zeta_1^{*}A \lto \zeta_2^{*}A
\end{equation*}
of vector bundles over $Z^{[2]} = Z \times_P Z$, where $P:=Y_1 \times_M Y_2$, with the
following properties:
\begin{itemize}
\item[a)]
It satisfies
the cocycle condition
\begin{equation*}
\zeta_{13}^{*}\mathrm{d}_\mathcal{A} = \zeta_{23}^{*}\mathrm{d}_\mathcal{A} \circ \zeta_{12}^{*}\mathrm{d}_\mathcal{A}
\end{equation*}
as an equation of isomorphisms of vector bundles over $Z^{[3]}$.
\item[b)]
The diagram
\begin{equation*}
\alxydim{@C=1.5cm@R=1.5cm}{L_1 \otimes \zeta_3^{*}A \ar[d]_{\id \otimes \zeta_{34}^{*}\mathrm{d}_{\mathcal{A}}}
\ar[r]^{\zeta_{13}^{*}\alpha} & \zeta_1^{*}A \otimes L_2 \ar[d]^{\zeta_{12}^{*}\mathrm{d}_{\mathcal{A}}
\otimes \id}
\\ L_1 \otimes \zeta_4^{*}A \ar[r]_{\zeta_{24}^{*}\alpha} & \zeta_{2}^{*}A
\otimes L_2}
\end{equation*} 
of isomorphisms of vector bundles over $Z^{[2]} \times_M Z^{[2]}$ is commutative.
\end{itemize}
\end{lemma} 

\proof
Notice that the isomorphism $\alpha$ of $\mathcal{A}$ restricted from $Z \times_M
Z$ to $Z \times_P Z$ gives an isomorphism
\begin{equation}
\alpha|_{Z\times_P Z}\lmaps \Delta^{*}L_1 \otimes \zeta_2^{*}A \lto \zeta_1^{*}A
\otimes \Delta^{*}L_2\text{.}
\end{equation}
By composition with the isomorphisms $t_{\mu_1}$ and $t_{\mu_2}$ from Lemma
\ref{lem4} we obtain the isomorphism $\mathrm{d}_\mathcal{A}$:
\begin{equation}
\alxydim{@C=1.7cm}{\zeta_1^{*}A \ar[r]^-{\id \otimes t_{\mu_2}^{-1}} & \zeta_1^{*} A\otimes
\Delta^{*}L_2 \ar[r]^-{\alpha|_{Z\times_P Z}^{-1}} & \Delta^{*}L_1 \otimes \zeta_2^{*}A \ar[r]^-{t_{\mu_1}
\otimes \id} & \zeta_2^{*}A\text{.}}
\end{equation} 
The cocycle condition a) and the commutative diagram b)  follow both from axiom \hyperlink{1M2}{(1M2)} for $\mathcal{A}$ and the properties of the isomorphisms $t_{\mu_1}$
and $t_{\mu_2}$ from Lemma
\ref{lem4}.
\endofproof

Now that we have defined the objects of $\mathfrak{Hom}(\mathcal{G}_1,\mathcal{G}_2)$,
we come to its morphisms.
 For two 1-morphisms $\mathcal{A}_1:\mathcal{G}_1 \to \mathcal{G}_2$
and $\mathcal{A}_2:\mathcal{G}_1 \to \mathcal{G}_2$,
consider  triples 
\begin{equation}
\label{9}
(W,\omega,\beta_W)
\end{equation}
 consisting of a smooth manifold $W$, a surjective submersion $\omega:W \to Z_1 \times_P Z_2$, where again $P:=Y_1 \times_M Y_2$, and
a morphism $\beta_W: A_1 \to A_2$ of vector bundles over
$W$.
Here we work again with the convention that we don't write pullbacks along
the unlabelled canonical projections $W \to Z_1$ and $W \to Z_2$. The triples (\ref{9}) have to satisfy one axiom \hypertarget{2M}{(2M)}: the isomorphism $\beta_W$ has to be compatible with isomorphism
$\alpha_1$ and $\alpha_2$ of the 1-morphisms $\mathcal{A}_1$ and $\mathcal{A}_2$
in the sense that the diagram 
\begin{equation}
\label{12}
\alxydim{@C=1.5cm@R=1.5cm}{
L_1 \otimes \omega_2^{*}A_1 \ar[r]^-{\alpha_1}
\ar[d]_{1 \otimes \omega_2^{*} \beta_W}
& \omega_1^{*}A_1  \otimes L_2
\ \ar[d]^{\omega_1^{*}\beta_W\otimes
1}  \\ L_1
\otimes \omega_2^{*}A_2 \ar[r]_-{\alpha_2}
& \omega_1^{*}A_2  \otimes L_2}
\end{equation}
of morphisms of vector bundles over $W \times_M W$ is commutative. 
On the set of all  triples (\ref{9}) satisfying this axiom we define an equivalence relation according to that
two triples
$(W,\omega,\beta_W)$ and $(W',\omega',\beta_{W'})$
are equivalent,
 if  there exists a smooth manifold $X$ with surjective submersions to $W$
 and $W'$ for which the diagram 
\begin{equation}
\label{13}
\alxydim{@C=0.1cm@R=0.8cm}{&X \ar[dl]\ar[dr]&\\W \ar[dr]_<<<{\omega}&&W'\ar[dl]^<<<{\omega'}\\&Z_1 \times_P Z_2&}
\end{equation}
of surjective submersions is commutative, and the  morphisms $\beta_W$ and $\beta_{W'}$ coincide when pulled back to $X$. 

\begin{definition}
\label{def2}
A 2-morphism $\beta: \mathcal{A}_1 \Rightarrow \mathcal{A}_2$ is an equivalence
class of triples $(W,\omega,\beta_W)$ satisfying axiom \hyperlink{2M
}{(2M)}.
\end{definition}

As an example of a 2-morphism we define the identity  2-morphism $\id_{\mathcal{A}}:\mathcal{A}
\Rightarrow \mathcal{A}$ associated to every 1-morphism $\mathcal{A}:\mathcal{G}_1
\to \mathcal{G}_2$. It is defined as the equivalence class of the triple $(Z^{[2]},\id_{Z^{[2]}},\mathrm{d}_{\mathcal{A}})$ consisting of the fibre product $Z^{[2]}=Z \times_P Z$, the
identity  $\id_{Z[2]}$ and the isomorphism $\mathrm{d}_\mathcal{A}:\zeta_1^{*}A
\to \zeta_2^{*}A$ of vector bundles over $Z^{[2]}$ from Lemma \ref{lem3}. Axiom \hyperlink{2M}{(2M)} for this
triple is proven with Lemma \ref{lem3} b). 

\medskip

Now we have defined objects and morphisms of the morphism category $\mathfrak{Hom}(\mathcal{G}_1,\mathcal{G}_2)$,
and we continue with the definition the composition $\beta' \bullet \beta$ of two 2-morphisms $\beta:\mathcal{A}_1
\Rightarrow \mathcal{A}_2$ and $\beta':\mathcal{A}_2 \Rightarrow \mathcal{A}_3$.
It is called  vertical composition in agreement with  the diagrammatical notation
\begin{equation}
\alxydim{@C=2cm}{\mathcal{G}_1 \ar@/^2.5pc/[r]^{\mathcal{A}_1}="1"
\ar[r]^<<<<<<<{\mathcal{A}_2} \ar@/_2.5pc/[r]_{\mathcal{A}_3}="3"
\ar@{}[r]^{}="2" & \mathcal{G}_2\ar@{=>}"1";"2"|{\beta}
\ar@{=>}"2";"3"|{\beta'}} \text{.}
\label{vercomp}
\end{equation}
We choose representatives $(W,\omega,\beta_W)$ and $(W',\omega',\beta_{W'})$
and consider the fibre product $\tilde W:= W\times_{Z_2}
W'$ with its canonical surjective submersion $\tilde \omega: \tilde W \to
Z_1 \times_P Z_3$, where again $P:=Y_1 \times_M Y_2$. By construction we can compose the pullbacks of the morphisms $\beta_W$
and $\beta_{W'}$ to $\tilde W$ and obtain a morphism
\begin{equation}
\beta_{W'} \circ \beta_W\lmaps A_1 \lto A_3
\end{equation}
of vector bundles over $\tilde W$. From  axiom \hyperlink{2M}{(2M)} for $\beta$ and $\beta'$
the one for the triple $(\tilde W,\tilde \omega, \beta_{W'}\circ \beta_W)$
follows. Furthermore, the equivalence class of this triple is independent
of the choices of the representatives of $\beta$ and $\beta'$and thus
defines  
 the 2-morphism  $\beta'\bullet\beta$. The composition
 $\bullet$ of the category $\mathfrak{Hom}(\mathcal{G}_1,\mathcal{G}_2)$
 defined like this is associative. 

It remains to check that the  2-isomorphism
 $\id_{\mathcal{A}}:\mathcal{A} \Rightarrow \mathcal{A}$ defined above is the identity under the composition $\bullet$. Let $\beta: \mathcal{A} \Rightarrow \mathcal{A}'$
be a 2-morphism and $(W,\omega,\beta_{W})$  a representative.  The composite $\beta \bullet \id_{\mathcal{A}}$ can be represented by the triple $(W',\omega',\beta \circ \mathrm{d}_{\mathcal{A}})$ with $W'=Z \times_P  W$, where $\omega':W' \to Z \times_P
Z'$ is the identity on the first factor and the projection $W \to Z'$ on the second one. We have to show, that this triple is equivalent to the original
representative $(W,\omega,\beta_W)$ of $\beta$. Consider the fibre product
\begin{equation}
X:= W \times_{(Z \times_P Z')} W' \cong W \times_{Z'}W\text{,}
\end{equation}
so that condition (\ref{13}) is satisfied. The restriction of
the commutative diagram (\ref{12}) of morphisms of vector bundles over
$W \times_M W$ from axiom \hyperlink{2M}{(2M)} for $\beta$ to $X$ gives rise
to the commutative diagram
\begin{equation}
\label{3}
\alxydim{@C=2cm@R=1cm}{
 \zeta_2^{*}A \ar[r]^-{\mathrm{d}_{\mathcal{A}}^{-1}}
\ar[d]_{ \omega_2^{*} \beta_W}
& \zeta_1^{*}A  
\ \ar[d]^{\omega_1^{*}\beta_W}  \\  A' \ar[r]_-{\Delta^{*}\mathrm{d}_{\mathcal{A}'}^{-1}}
& A'  }
\end{equation}
of morphisms of vector bundles over $X$, where $\mathrm{d}_{\mathcal{A}}$ and $\mathrm{d}_{\mathcal{A}'}$
are the isomorphisms determined by the 1-morphisms $\mathcal{A}$ and $\mathcal{A}'$
according to Lemma \ref{lem3}. Their cocycle condition from Lemma \ref{lem3} a)  implies $\Delta^{*}\mathrm{d}_{\mathcal{A}'}=\id$,
so that  diagram (\ref{3}) is reduced to the equality $\omega_2^{*}\beta_W \circ \mathrm{d}_{\mathcal{A}}
= \omega_1^{*}\beta_W$ of isomorphisms of vector bundles over $X$. This shows
that  the triples $(W,\omega,\beta_W)$
and $(W',\omega',\beta_W \circ \mathrm{d}_{\mathcal{A}})$ are equivalent and we have $\beta
\bullet \id_{\mathcal{A}}=\beta$. The equality $\id_{\mathcal{A}'} \bullet
\beta= \beta$ follows analogously.

\medskip

Now the definition of the morphism category $\mathfrak{Hom}(\mathcal{G}_1,\mathcal{G}_2)$
is complete. A morphism in this category, i.e. a 2-morphism $\beta:\mathcal{A}\Rightarrow \mathcal{A}'$, is invertible if and only if the morphism $\beta_W: A \to
A'$ of any representative $(W,\omega,\beta_W)$ of $\beta$ is invertible.
Since
-- following our convention -- 
morphism of vector bundles respect the hermitian structures, $\beta_W$ is
invertible if and only if  the ranks of the vector bundles of the 1-morphisms
$\mathcal{A}$ and $\mathcal{A}'$ coincide.
In the following, we call two 1-morphisms $\mathcal{A}:\mathcal{G}_1 \to \mathcal{G}_2$
and $\mathcal{A}':\mathcal{G}_1 \to \mathcal{G}_2$ isomorphic, if there exists
a 2-isomorphism $\beta:\mathcal{A}\Rightarrow \mathcal{A}'$ between them.

\subsection{The Composition Functor}

\label{sec2}

Let $\mathcal{G}_1$, $\mathcal{G}_2$ and $\mathcal{G}_3$ be three bundles
gerbes over $M$. We define the composition functor
\begin{equation}
\circ \lmaps \mathfrak{Hom}(\mathcal{G}_2,\mathcal{G}_3) \times \mathfrak{Hom}(\mathcal{G}_1,\mathcal{G}_2)
\longrightarrow \mathfrak{Hom}(\mathcal{G}_1,\mathcal{G}_3)
\end{equation}
on objects in the following way.
Let $\mathcal{A}:\mathcal{G}_1
\to \mathcal{G}_2$ and $\mathcal{A}':\mathcal{G}_2
\to \mathcal{G}_3$ be two 1-morphisms. The composed 1-morphism 
\begin{equation}
\mathcal{A}'
\circ \mathcal{A}\lmaps\mathcal{G}_1
\lto \mathcal{G}_3
\end{equation}
consists of the fibre
product $\tilde Z := Z \times_{Y_2}
Z'$ with its canonical surjective submersion $\tilde \zeta: \tilde Z \to Y_1
\times_M Y_3$, 
 the
vector bundle $\tilde A := A \otimes A'$ over $\tilde Z$, and the isomorphism
\begin{equation}
\tilde \alpha := (\id_{\zeta_1^{*}A} \otimes
\alpha') \circ
(\alpha \otimes \id_{\zeta_2^{\prime*}A'})
\end{equation}
of vector bundles over $\tilde
Z \times_M \tilde Z$. 

Indeed, this defines a 1-morphism from $\mathcal{G}_1$ to $\mathcal{G}_3$.
Recall that if $\nabla_A$ and $\nabla_{A'}$ denote the connections on the vector bundles
$A$ and $A'$, the tensor product connection $\nabla$ on $A \otimes A'$ is
defined by
\begin{equation}
\nabla(\sigma \otimes \sigma') = \nabla_A(\sigma) \otimes \sigma' + \sigma
\otimes \nabla_{A'}(\sigma')
\end{equation}
for sections $\sigma\in\Gamma(A)$ and $\sigma'\in\Gamma(A')$.
If we take $n$ to be the rank of $A$ and $n'$ the rank of $A'$ the curvature of the tensor product vector bundle is 
\begin{equation}
\mathrm{curv}(A \otimes A') = \mathrm{curv}(A) \otimes \id_{n'} + \id_{n}
\otimes \mathrm{curv}(A')\text{.}
\end{equation}
Hence its trace
\begin{eqnarray}
\label{36}
\frac{1}{nn'}\mathrm{tr}(\mathrm{curv}(\tilde A)) &=&\frac{1}{n} \mathrm{tr}(\mathrm{curv}(A))
+\frac{1}{n'}\mathrm{tr}( \mathrm{curv}(A'))\nonumber\\ &=& C_2 - C_1
+ C_3 - C_2 \nonumber \\ &=& C_3 - C_1
\end{eqnarray}
satisfies axiom \hyperlink{1M1}{(1M1)}. Notice that equation (\ref{36}) involves
unlabeled projections from $\tilde Z$ to $Y_1$, $Y_2$ and $Y_3$, where the
one to $Y_2$ is  unique because $\tilde Z$ is the fibre product
over $Y_2$.  Furthermore,  $\tilde
\alpha$ is  an isomorphism 
\begin{equation}
\alxy{L_1 \otimes \tilde \zeta_2^{*}\tilde A \ar@{=}[r] & L_1 \otimes \zeta_2^{*}A \otimes
\zeta_2^{\prime*}A' \ar[d]^{\alpha
\otimes \id} & \\ &
\zeta_1^{*}A \otimes L_2 \otimes
\zeta_2^{\prime*}A' \ar[d]^{\id
\otimes \alpha'} & \\ &
\zeta_1^{*}A \otimes \zeta_1^{\prime*}A'
\otimes L_3 \ar@{=}[r] & \tilde
\zeta_1^{*}\tilde A \otimes L_3\text{.}}
\end{equation}
Axiom \hyperlink{1M2}{(1M2)} follows from axioms \hyperlink{1M2}{(1M2)}
for $\mathcal{A}$ and $\mathcal{A}'$.

\begin{proposition}
\label{lem2}
The composition of 1-morphisms 
is strictly associative: for three 1-morphisms $\mathcal{A}:\mathcal{G}_1 \to
\mathcal{G}_2$, $\mathcal{A}':\mathcal{G}_2 \to \mathcal{G}_3$ and $\mathcal{A}'':\mathcal{G}_3
\to \mathcal{G}_4$ we have
\begin{equation*}
(\mathcal{A}'' \circ \mathcal{A}') \circ \mathcal{A} = \mathcal{A}'' \circ
(\mathcal{A}' \circ \mathcal{A})\text{.}
\end{equation*}
\end{proposition}

\proof
By definition, both 1-morphism $(\mathcal{A}'' \circ \mathcal{A}') \circ \mathcal{A} $ and $\mathcal{A}'' \circ
(\mathcal{A}' \circ \mathcal{A})$ consist of the smooth manifold $X=Z \times_{Y_2} Z' \times_{Y_3} Z''$ with the same surjective submersion $X \to Y_1 \times_M
Y_4$. On $X$,
they have the same vector bundle $A \otimes A' \otimes A''$, and finally the
same isomorphism
\begin{equation}
(\id \otimes \id \otimes \alpha'') \circ (\id \otimes \alpha' \otimes \id) \circ (\alpha \otimes \id \otimes \id)
\end{equation}
of vector bundles over $X \times_M X$.
\endofproof

Now we have to define the functor $\circ$ on 2-morphisms.
Let $\mathcal{A}_1,\mathcal{A}_1':
\mathcal{G}_1 \to \mathcal{G}_2$ and $\mathcal{A}_2,\mathcal{A}_2':\mathcal{G}_2
\to \mathcal{G}_3$ be 1-morphisms between bundle gerbes. 
The functor $\circ$ on morphisms is called horizontal composition due to the diagrammatical notation
\begin{equation}
\alxydim{@C=1cm}{\mathcal{G}_1 \ar@/^1.6pc/[rr]^{\mathcal{A}_1}="1"  \ar@/_1.6pc/[rr]_{\mathcal{A}'_1}="2"
 \ar@{=>}"1";"2"|{\beta_1} && \mathcal{G}_2 \ar@/^1.6pc/[rr]^{\mathcal{A}_2}="3"  \ar@/_1.6pc/[rr]_{\mathcal{A}'_2}="4"
 \ar@{=>}"3";"4"|{\beta_2}
&& \mathcal{G}_3} \;\;=\;\;  \alxydim{}{\mathcal{G}_1
\ar@/^1.6pc/[rrr]^{\mathcal{A}_2
\circ \mathcal{A}_1}="5"  \ar@/_1.6pc/[rrr]_{\mathcal{A}'_2
\circ \mathcal{A}'_1}="6"
 \ar@{=>}"5";"6"|{\beta_2 \circ \beta_1}
&&& \mathcal{G}_3}\text{.}
\label{horcomp}
\end{equation}
Recall that the compositions $\mathcal{A}_2 \circ
\mathcal{A}_1$ and $\mathcal{A}_2' \circ \mathcal{A}_1'$ consist of smooth
manifolds $\tilde Z = Z_1 \times_{Y_2} Z_2$ and $\tilde Z' = Z_1' \times
_{Y_2} Z_2'$ with surjective submersions to $P:=Y_1 \times_M Y_3$,  of vector bundles $\tilde A := A_1 \otimes A_2$
over $\tilde Z$ and $\tilde A':=A_1' \otimes A_2'$ over $\tilde Z'$, and
of isomorphisms $\tilde \alpha$ and $\tilde \alpha'$ over $\tilde Z \times_M
\tilde Z$
and $\tilde Z' \times_M \tilde Z'$.

To define
the composed 2-morphism $\beta_2 \circ \beta_1$, we first need a surjective submersion
\begin{equation}
\omega\lmaps W \lto \tilde Z \times_P \tilde Z'\text{.}
\end{equation}
We choose representatives 
$(W_1,\omega_1,\beta_{W_{1}})$ and $(W_2,\omega_2,\beta_{W_{2}})$ of the
2-morphisms $\beta_1$ and $\beta_2$ and define
\begin{equation}
\label{20}
W:=  \tilde Z \times_{P} (W_1 \times_{Y_2} W_2) \times_P \tilde Z'
\end{equation}
with the surjective submersion $\omega:= \tilde z \times \tilde z'$ projecting on the  first and the last factor. Then, we need a morphism $\beta_{W}: \tilde z^{*}\tilde A \to
\tilde z'^{*} \tilde A'$ of vector bundles over $W$. Notice that we have  maps
\begin{equation}
u\lmaps W_1 \times_{Y_2} W_2 \lto \tilde Z
\quad\text{ and }\quad
u'\lmaps W_1 \times_{Y_2} W_2 \lto \tilde Z'
\end{equation}
such that we obtain surjective submersions  
\begin{equation}
\tilde z \times u\lmaps W \lto \tilde Z^{[2]}
\quad\text{ and }\quad
u' \times \tilde z'\lmaps W \lto \tilde Z'^{[2]}.
\end{equation}
Recall from Lemma \ref{lem3} that the 1-morphisms $\mathcal{A}_2 \circ \mathcal{A}_1$ and $\mathcal{A}_2' \circ \mathcal{A}_1'$ define
 isomorphisms $\mathrm{d}_{\mathcal{A}_2\circ \mathcal{A}_1}$ and $\mathrm{d}_{\mathcal{A}'_2\circ \mathcal{A}'_1}$ of vector bundles over $\tilde Z^{[2]}$ and $\tilde Z'^{[2]}$,
whose pullbacks
to $W$ along the above maps
 are
isomorphisms\begin{equation}
\mathrm{d}_{\mathcal{A}_2\circ \mathcal{A}_1}\lmaps \tilde z^{*}\tilde A \lto  u^{*}\tilde A
\quad\text{ and }\quad
\mathrm{d}_{\mathcal{A}'_2\circ \mathcal{A}'_1}\lmaps u'^{*}\tilde A' \lto \tilde z'^{*}\tilde A'
\end{equation}
of vector bundles over $W$. 
Finally, the morphisms $\beta_{W_1}$ and $\beta_{W_2}$ give a morphism
\begin{equation}
\tilde\beta := \beta_{W_1} \otimes \beta_{W_2}\lmaps u^{*}\tilde A \lto u'^{*}\tilde A'
\end{equation}
of vector bundles over $W$ so that the composition
\begin{equation}
\label{22}
\beta_{W}:=\mathrm{d}_{\mathcal{A}'_2\circ \mathcal{A}'_1} \circ \tilde \beta \circ \mathrm{d}_{\mathcal{A}_2\circ \mathcal{A}_1}
\end{equation}
is a well-defined morphism of vector bundles over $W$. Axiom \hyperlink{2M}{(2M)} for the triple $(W,\omega,\beta_W)$ follows from Lemma \ref{lem3} b) for  $\mathcal{A}_2 \circ \mathcal{A}_1$ and $\mathcal{A}_2' \circ \mathcal{A}_1'$ and from the axioms \hyperlink{2M}{(2M)} for the representatives of $\beta_1$
and $\beta_2$. Furthermore, the equivalence class of $(W,\omega,\beta_W)$ is
independent of the choices of the representatives of $\beta_1$ and $\beta_2$.

\begin{lemma}
\label{lem5}
The assignment $\circ$, defined above on objects and morphisms, is a functor
\begin{equation*}
\circ \lmaps \mathfrak{Hom}(\mathcal{G}_2,\mathcal{G}_3) \times \mathfrak{Hom}(\mathcal{G}_1,\mathcal{G}_2)
\longrightarrow \mathfrak{Hom}(\mathcal{G}_1,\mathcal{G}_3)\text{.}
\end{equation*}
\end{lemma}
  
\proof
i) The assignment $\circ$ respects identities, i.e. for 1-morphisms $\mathcal{A}_1:\mathcal{G}_1
\to \mathcal{G}_2$ and $\mathcal{A}_2:\mathcal{G}_2 \to \mathcal{G}_3$,
\begin{equation}
\id_{\mathcal{A}_2} \circ \id_{\mathcal{A}_1} = \id_{\mathcal{A}_2 \circ
\mathcal{A}_1}\text{.}
\end{equation}
To show this we choose the defining representatives $(W_1,\id,\mathrm{d}_{\alpha_1})$ of $\id_{\mathcal{A}_1}$ and $(W_2,\id,\mathrm{d}_{\alpha_2})$ of $\id_{\mathcal{A}_2}$, where $W_1= Z_1 \times_{(Y_1 \times_M
Y_2)} Z_1$ and $W_2=Z_2 \times_{(Y_2 \times_M Y_3)} Z_2$.
Consider the diffeomorphism
\begin{equation}
\label{21}
f:W_1
\times_{Y_2} W_2 \to \tilde Z \times_{Y_1 \times_M Y_2 \times_M Y_3} \tilde
Z: (z_1,z_1',z_2,z_2') \mapsto (z_1,z_2,z_1',z_2')\text{,}
\end{equation}
where $\tilde Z= Z_1 \times_{Y_2} Z_2$. 
From the definitions of the isomorphisms $\mathrm{d}_{\mathcal{A}_1}$, $\mathrm{d}_{\mathcal{A}_2}$
 and $\mathrm{d}_{\mathcal{A}_2 \circ \mathcal{A}_2}$ we conclude the
 equation
\begin{equation}
\label{25}
\mathrm{d}_{\mathcal{A}_{1}} \otimes \mathrm{d}_{\mathcal{A}_{2}} = f^{*}\mathrm{d}_{\mathcal{A}_2\circ \mathcal{A}_1}
\end{equation}
of isomorphisms of vector bundles over $W_1 \times_{Y_2} W_2$.
The horizontal composition $\id_{\mathcal{A}_2}
\circ \id_{\mathcal{A}_1}$ is canonically represented by the triple $(W,\omega,\beta_W)$
where $W$ is defined in (\ref{20}) and $\beta_W$ is defined in (\ref{22}).
Now, the diffeomorphism $f$ extends to an embedding $f:W \to \tilde Z^{[4]}$
 into the four-fold
fibre product of $\tilde Z$ over $P=Y_1 \times_M Y_3$, such that $\omega:W
\to \tilde Z^{[2]}$ factorizes
over $f$,
\begin{equation}
\label{27}
\omega=\tilde\zeta_{14}\circ f\text{.}
\end{equation}
From (\ref{22}) and (\ref{25}) we obtain
\begin{eqnarray}
\beta_W  &=&
 \mathrm{d}_{\mathcal{A}_2\circ \mathcal{A}_1} \circ (\mathrm{d}_{\mathcal{A}_{1}} \otimes \mathrm{d}_{\mathcal{A}_{2}})
\circ \mathrm{d}_{\mathcal{A}_2\circ \mathcal{A}_1} \nonumber \\&=& f^{*}(\tilde\zeta_{34}^{*}\mathrm{d}_{\mathcal{A}_2\circ \mathcal{A}_1} \circ \tilde\zeta_{23}^{*}\mathrm{d}_{\mathcal{A}_2\circ \mathcal{A}_1}
\circ \tilde\zeta_{12}^{*}\mathrm{d}_{\mathcal{A}_2\circ \mathcal{A}_1}).
\end{eqnarray}
The cocycle condition for $\mathrm{d}_{\mathcal{A}_2\circ \mathcal{A}_1}$ from
Lemma \ref{lem3} a) and  (\ref{27}) give 
\begin{equation}
\label{26}
\beta_W=f^{*}\tilde\zeta_{14}^{*}\mathrm{d}_{\mathcal{A}_2\circ \mathcal{A}_1}=\omega^{*}\mathrm{d}_{\mathcal{A}_2\circ \mathcal{A}_1}\text{.}
\end{equation}
We had to show that the triple $(W,\omega,\beta_W)$ which represents $\id_{\mathcal{A}_2
}\circ \id_{\mathcal{A}_2}$ is equivalent to the
triple $(\tilde Z^{[2]},\id,\mathrm{d}_{\mathcal{A}_2\circ \mathcal{A}_1})$ which
defines the identity 2-morphism $\id_{\mathcal{A}_2 \circ \mathcal{A}_1}$.
For the choice $X:= W$ with surjective submersions $\id:X \to W$ and $\omega:X
\to \tilde Z^{[2]}$, equation (\ref{26}) shows exactly this equivalence.

ii) The assignment $\circ$ respects the composition $\bullet$, i.e. for 2-morphisms
$\beta_i: \mathcal{A}_i \Rightarrow \mathcal{A}_i'$ and $\beta_i':\mathcal{A}_i'
\Rightarrow \mathcal{A}_i''$ between 1-morphisms $\mathcal{A}_i$, $\mathcal{A}_i'$
and $\mathcal{A}_i''$ from   $\mathcal{G}_i$ to $\mathcal{G}_{i+1}$, everything for $i=1,2$, we have an equality
\begin{equation}
\label{29}
(\beta_2' \bullet \beta_2) \circ (\beta_1' \bullet \beta_1) = (\beta_2' \circ
\beta_1') \bullet (\beta_2 \circ \beta_1)
\end{equation}
of 2-morphisms from $\mathcal{A}_2 \circ \mathcal{A}_1$ to $\mathcal{A}_2''
\circ \mathcal{A}_1''$. 
This equality is also known as the compatibility of vertical and horizontal
compositions. To prove it, let us introduce the notation $\tilde Z:= Z_1 \times_{Y_2}
Z_2$, and analogously $\tilde Z'$ and $\tilde Z''$, furthermore we write $P:=Y_1
\times_M Y_3$. 
Notice that the 2-morphism on the left hand side
of (\ref{29}) is represented by a triple $(V,\nu,\beta_V)$ with 
\begin{equation}
V= \tilde Z \times _{P} (\tilde W_1 \times_{Y_2} \tilde W_2) \times_{P}
\tilde Z''\text{,}
\end{equation}
where the fibre products $\tilde W_i := W_i \times_{Z_i'} W_i'$ arise from
the vertical compositions $\beta_i' \bullet \beta_i$. The surjective submersion $\nu:V \to \tilde Z
\times_{P} \tilde Z''$ is the projection on the first and the last factor, and
\begin{equation}
\beta_V = \mathrm{d}_{\mathcal{A}_2'' \circ \mathcal{A}_1''} \circ ((\beta_1'
\circ \beta_1) \otimes
(\beta_2' \circ \beta_2)) \circ \mathrm{d}_{\mathcal{A}_2 \circ \mathcal{A}_1}
\end{equation}
is a morphism of vector bundles over $V$. The 2-morphism on the right
hand side of (\ref{29}) is represented by the triple $(V',\nu',\beta_{V'})$
with
\begin{eqnarray}
V' &=& (\tilde Z \times_{P} (W_1 \times_{Y_2} W_2) \times_P \tilde Z') \times_{\tilde
Z'} (\tilde Z' \times_P (W'_1 \times_{Y_2} W_2') \times_P \tilde Z'') \nonumber
\\&\cong& \tilde Z \times_{P} (W_1 \times_{Y_2} W_2) \times_P \tilde Z' \times_P (W'_1 \times_{Y_2} W_2') \times_P \tilde Z''\text{,}
\end{eqnarray}
where $\nu'$ is again the projection on the outer factors, and 
\begin{equation}
\beta_{V'} = \mathrm{d}_{\mathcal{A}_2'' \circ \mathcal{A}_1''} \circ (\beta_1'
\otimes \beta_2') \circ \mathrm{d}_{\mathcal{A}_2' \circ \mathcal{A}_1'}
\circ (\beta_1 \otimes \beta_2) \circ \mathrm{d}_{\mathcal{A}_2 \circ \mathcal{A}_1}\text{,}
\end{equation}
where we have used the cocycle condition for $\mathrm{d}_{\mathcal{A}_2' \circ
\mathcal{A}_1'}$ from Lemma \ref{lem3} b). 

We have to show that the triples $(V,\nu,\beta_{V})$ and $(V',\nu',\beta_{V'})$
are equivalent. Consider the fibre product
\begin{equation}
X:= V \times_{\tilde Z \times_P \tilde Z''} V'
\end{equation}
with surjective submersions $v:X \to V$ and $v':X \to V'$. To show the equivalence
of the two triples, we have to prove the equality
\begin{equation}
v^{*}\beta_V = v'^{*}\beta_{V'}\text{.}
\end{equation}
It is equivalent to the commutativity of the outer shape of the following
diagram of isomorphisms of vector bundles over $X$: 
\begin{equation}
\alxydim{@C=0.7cm@R=0.5cm}{&A_1 \otimes A_2 \ar[dl]_{\mathrm{d}_{\mathcal{A}_2 \circ \mathcal{A}_1}}
\ar[dr]^{\mathrm{d}_{\mathcal{A}_2 \circ \mathcal{A}_1}} & \\ v^{*}(A_1 \otimes
A_2) \ar[rr]|{\mathrm{d}_{\mathcal{A}_2 \circ \mathcal{A}_1}} \ar[dd]_{\beta_1
\otimes \beta_2} && v'^{*}(A_1 \otimes A_2) \ar[d]^{\beta_1 \otimes \beta_2} \\ &&
v'^{*}(A_1' \otimes A_2') \ar[dd]^{\mathrm{d}_{\mathcal{A}_2' \circ \mathcal{A}_1'}} \\
v^{*}(A_1'
\otimes A_2') \ar[dd]_{\beta_1' \otimes \beta_2'} \ar[rru]|{\mathrm{d}_{\mathcal{A}_2'\circ \mathcal{A}_1'}} \ar[rrd]|{\mathrm{d}_{\mathcal{A}_2'
\circ \mathcal{A}_1'}} && \\ && v'^{*}(A_1' \otimes A_2') \ar[d]^{\beta_1' \otimes
\beta_2'} \\ v^{*}(A_1'' \otimes A_2'') \ar[rr]|{\mathrm{d}_{\mathcal{A}_2'' \circ
\mathcal{A}_1''}} \ar[dr]_{\mathrm{d}_{\mathcal{A}_2'' \circ \mathcal{A}_1''}} && v'^{*}(A_1'' \otimes A_2'') \ar[dl]^{\mathrm{d}_{\mathcal{A}_2'' \circ \mathcal{A}_1''}}
\\ &A_1'' \otimes A_2''&}
\end{equation}
The commutativity of the outer shape of this diagram follows from the commutativity
of its five subdiagrams: the triangular ones are commutative due to the cocycle
condition from Lemma \ref{lem3} a), and the commutativity of the foursquare ones follows from  axiom \hyperlink{2M}{(2M)} of the 2-morphisms. 
\endofproof

\medskip

To finish the definition of the 2-category $\mathfrak{BGrb}(M)$ we have to define the natural 2-isomorphisms $\lambda_{\mathcal{A}}: \mathcal{A}\circ \id_{\mathcal{G}}
\Rightarrow \mathcal{A}$ and $\rho_{\mathcal{A}}: \id_{\mathcal{G}'} \circ
\mathcal{A} \Rightarrow \mathcal{A}$ for a given 1-morphism $\mathcal{A}:\mathcal{G}
\to \mathcal{G}'$, and we have to show that they satisfy axiom \hyperlink{2C2}{(2C2)}. We define the 1-morphism $\mathcal{A}\circ \id_{\mathcal{G}}$
as follows: it has the canonical surjective submersion from $\tilde Z = Y^{[2]} \times_{Y} Z \cong Y \times_M
Z$ to $P:=Y \times_M Y'$ and the vector bundle $L \otimes A$ over $\tilde Z$. Consider 
\begin{equation}
W := \tilde Z \times_P Z \cong Z \times_{Y'}Z
\end{equation}
and the identity
$\omega := \id_{W}$. Under this identification, let us consider the restriction of the isomorphism $\alpha$ of the 1-morphism
$\mathcal{A}$ from $Z \times_M Z$ to $W=Z \times_{Y'} Z$. If $s:W \to W$
denotes the exchange of the two factors, we obtain an isomorphism
\begin{equation}
s^{*}\alpha|_{W}\lmaps L \otimes \zeta_1^{*}A \lto \zeta_2^{*}A \otimes \Delta^{*}L'
\end{equation}
of vector bundles over $W$. By composition with the canonical trivialization
of the line bundle $\Delta^{*}L'$ from Lemma \ref{lem4} it gives an isomorphism \begin{equation}
\lambda_W := (\id \otimes t_{\mu'}) \circ s^{*}\alpha|_W\lmaps  L \otimes
\zeta_1^{*}A \lto \zeta_2^{*}A 
\end{equation}
of vector bundles over $W$. The axiom \hyperlink{2M}{(2M)} for the triple $(W,\omega,\lambda_W)$
follows from axiom \hyperlink{1M2}{(1M2)} for
the 1-morphism $\mathcal{A}$ and from the properties of $t_{\mu'}$ from Lemma
\ref{lem4}. So, $\lambda_{\mathcal{A}}$ is defined to be the equivalence
class of this triple. The definition of $\rho_{\mathcal{A}}$ goes analogously:
we take $W=Z \times_Y Z$ and obtain by restriction the isomorphism
\begin{equation}
\alpha|_W\lmaps \Delta^{*}L \otimes \zeta_2^{*}A \lto \zeta_1^{*}A \otimes L'. 
\end{equation}
Then, the 2-isomorphism $\rho_{\mathcal{A}}$ is defined by the triple $(W,\omega,\rho_W)$
with the isomorphism
\begin{equation}
\rho_W := (t_{\mu} \otimes \id) \circ \alpha|_W^{-1}\lmaps \zeta_1^{*}A  \otimes
L '\lto 
\zeta_2^{*}A 
\end{equation}
of vector bundles over $W$. 

\begin{lemma}
The 2 -isomorphisms $\lambda_{\mathcal{A}}$
and $\rho_{\mathcal{A}}$ are natural in $\mathcal{A}$, i.e. for any 2-morphism $\beta:\mathcal{A}\Rightarrow
\mathcal{A}'$ the naturality squares
\begin{equation*}
\alxydim{@C=1.5cm@R=1cm}{\id_{\mathcal{G}'} \circ \mathcal{A} \ar@{=>}[d]_{\id_{\id_{\mathcal{G}'}}
\circ \beta} \ar@{=>}[r]^-{\rho_{\mathcal{A}}} & \mathcal{A} \ar@{=>}[d]^{\beta}
\\ \id_{\mathcal{G}'} \circ \mathcal{A}' \ar@{=>}[r]_-{\rho_{\mathcal{A}'}}
& \mathcal{A}'}
\quad\text{ and }\quad
\alxydim{@C=1.5cm@R=1cm}{\mathcal{A} \circ \id_{\mathcal{G}}  \ar@{=>}[d]_{\beta
\circ \id_{\id_{\mathcal{G}}}} \ar@{=>}[r]^-{\lambda_{\mathcal{A}}} & \mathcal{A} \ar@{=>}[d]^{\beta}
\\ \mathcal{A}' \circ \id_{\mathcal{G}} \ar@{=>}[r]_-{\lambda_{\mathcal{A}'}}
& \mathcal{A}'}
\end{equation*}
are commutative. 
\end{lemma}

\proof
To calculate for instance the horizontal composition
$\id_{\id_{\mathcal{G}'}}\circ
\beta$ in the diagram on the left hand side first note that $\id_{\id_{\mathcal{G}'}}$
is canonically represented by the triple $(Y'^{[2]},\id,\id_L)$. The isomorphism
\begin{equation}
\mathrm{d}_{\id_{\mathcal{G}'} \circ \mathcal{A}}: \tilde\zeta_1^{*}(A
\otimes L') \to \tilde\zeta_2^{*}(A \otimes L')\text{,}
\end{equation}
which appears in the definition of
the horizontal composition, is an isomorphism of vector bundles over $\tilde Z \times_{Y \times_M Y'} \tilde Z$, where $\tilde\zeta:\tilde Z:=Z \times_M
Y' \to Y \times_M Y'$ is the surjective submersion of the composite
$\id_{\mathcal{G}'} \circ \mathcal{A}$. Here it simplifies to
\begin{equation}
\mathrm{d}_{\id_{\mathcal{G}'}
\circ \mathcal{A}}=(t_{\mu} \otimes \id \otimes \id) \circ (\alpha^{-1} \otimes
\id) \circ (1 \otimes \tilde\zeta_1^{*}\mu'^{-1})\text{.}
\end{equation}
With these simplifications and with axiom \hyperlink{1M2}{(1M2)} for $\mathcal{A}$
and $\mathcal{A}'$, the naturality squares  reduce to the compatibility
axiom \hyperlink{2M}{(2M)} of $\beta$ with the isomorphisms $\alpha$ and
$\alpha'$ of $\mathcal{A}$ and $\mathcal{A}'$ respectively. 
\endofproof

It remains to show that the  isomorphisms $\lambda_{\mathcal{A}}$
and $\rho_{\mathcal{A}}$ satisfy axiom \hyperlink{2C2}{(2C2)} of a 2-category.

\begin{proposition}
\label{lem1}
For 1-morphisms $\mathcal{A}:\mathcal{G}_1 \to \mathcal{G}_2$ and $\mathcal{A}':\mathcal{G}_2
\to \mathcal{G}_3$, the 2-isomorphisms $\lambda_{\mathcal{A}}$ and $\rho_{\mathcal{A}}$ satisfy
\begin{equation*}
\id_{\mathcal{A}'} \circ \rho_{\mathcal{A}} = \lambda_{\mathcal{A}'} \circ
\id_{\mathcal{A}}\text{.}
\end{equation*}
\end{proposition}

\proof 
The equation to prove is an equation of 2-morphisms from $\mathcal{A}' \circ
\id_{\mathcal{G}_2} \circ \mathcal{A}$ to $\mathcal{A}'\circ \mathcal{A}$.
The first 1-morphism consists of the surjective submersion $\tilde
Z := Z \times_M Z' \to P_{13}$, where we define $P_{ij}:=Y_i \times_M Y_j$,
further
of the vector bundle $A \otimes L_2 \otimes A'$ over $\tilde Z$. The second
1-morphism $\mathcal{A}'\circ \mathcal{A}$ consists of the surjective submersion
$\tilde Z':= Z \times_{Y_2} Z' \to P_{13}$ and the vector bundle $A \otimes
A'$ over $\tilde Z'$. Let us choose the defining representatives for the
involved
2-morphisms: we choose $(Z'^{[2]},\id,\mathrm{d}_{\mathcal{A}'})$ for $\id_{\mathcal{A}'}$,
with $W:= Z \times_{Y_1} Z$ we choose $(W,\id,\rho_W)$ for $\rho_{\mathcal{A}}$,
with $W':= Z' \times_{Y_3} Z'$ we choose $(W',\id,\lambda_{W'})$ for $\lambda_{\mathcal{A}'}$,
and we choose $(Z^{[2]},\id,\mathrm{d}_{\mathcal{A}'})$
for $\id_{\mathcal{A}}$. 

Now, the horizontal composition $\id_{\mathcal{A}'} \circ \rho_{\mathcal{A}}$
is defined by the triple $(V,\nu,\beta_V)$ with 
\begin{equation}
V=\tilde Z \times_{P_{13}} (W \times_{Y_2} Z'^{[2]}) \times_{P_{13}} \tilde
Z'\text{,}
\end{equation}
the projection $\nu:V \to \tilde Z \times_{P_{13}} \tilde Z'$ on the first and the last factor, and the isomorphism 
\begin{equation}
\beta_{V} = \mathrm{d}_{\mathcal{A}' \circ \mathcal{A}} \circ (\rho_W \otimes \mathrm{d}_{\mathcal{A}'}) \circ \mathrm{d}_{\mathcal{A}' \circ \id \circ \mathcal{A}}
\end{equation}
of vector bundles over $V$. The horizontal composition $\lambda_{\mathcal{A}'} \circ \id_{\mathcal{A}}$
is defined by the triple $(V',\nu',\beta_{V'})$ with
\begin{equation}
V' = \tilde Z \times_{P_{13}} (Z^{[2]} \times_{Y_2} W') \times_{P_{13}} \tilde
Z'\text{,} 
\end{equation}
again the projection $\nu'$ on the first and the last factor, and the isomorphism
\begin{equation}
\beta_{V'} = \mathrm{d}_{\mathcal{A}' \circ \mathcal{A}} \circ (\mathrm{d}_{\mathcal{A}}
\otimes \lambda_{W'}) \circ \mathrm{d}_{\mathcal{A}' \circ \id \circ \mathcal{A}}
\end{equation}
of vector bundles over $V$. 

To prove the proposition, we show that the triples $(V,\nu,\beta_V)$
and $(V',\nu',\beta_{V'})$ are equivalent. Consider the fibre product 
\begin{equation}
X
:= V \times_{(\tilde Z \times_{P_{13}} \tilde Z')} V'
\end{equation}
with surjective submersions $v:X \to V$ and $v':X \to V'$. 
The equivalence of the two triples follows from the equation
\begin{equation}
v^{*}\beta_V = v'^{*}\beta'
\end{equation}
of isomorphisms of vector bundles over $X$. It is equivalent to the commutativity of the outer
shape of the  following diagram of isomorphisms of vector bundles over $X$:
\begin{equation}
\alxydim{@R=1.2cm@C=0.2cm}{&A \otimes L_2 \otimes A' \ar[dl]_{\mathrm{d}_{\mathcal{A}'\circ \id
\circ \mathcal{A}}} \ar[dr]^{\mathrm{d}_{\mathcal{A}'\circ \id
\circ \mathcal{A}}}&\\ v^{*}(A \otimes L_2 \otimes A') \ar[d]_{\rho_{W}
\otimes \mathrm{d}_{\mathcal{A}'}} \ar[rr]|{\mathrm{d}_{\mathcal{A}'\circ \id
\circ \mathcal{A}}} && v'^{*}(A \otimes L_2 \otimes A') \ar[d]^{\mathrm{d}_{\mathcal{A}} \otimes \lambda_{W'}}
\\ v^{*}(A \otimes A') \ar[dr]_{\mathrm{d}_{\mathcal{A}' \circ \mathcal{A}}} \ar[rr]|{\mathrm{d}_{\mathcal{A}'\circ \mathcal{A}}} && v'^{*}(A \otimes A') \ar[dl]^{\mathrm{d}_{\mathcal{A}' \circ \mathcal{A}}} \\ &A \otimes A'&}
\end{equation}
The diagram is patched together from three subdiagrams, and the commutativity
of the outer shape follows because the three subdiagrams are commutative:
the triangle diagrams are commutative due to the cocycle condition from Lemma
\ref{lem3} b) for the 1-morphisms $\mathcal{A}' \circ \id_{\mathcal{G}_2}
\circ \mathcal{A}$ and $\mathcal{A}' \circ \mathcal{A}$ respectively. The
commutativity of the rectangular diagram in the middle follows from Lemma
\ref{lem4} and from axioms \hyperlink{1M2}{(1M2)} for $\mathcal{A}$ and $\mathcal{A}'$.
\endofproof

\subsection{Invertible 1-Morphisms}

\label{sec6}

In this subsection we address the question, which of the 1-morphisms of the
2-category $\mathfrak{BGrb}(M)$ are invertible. 
Let $\mathcal{G}_1$ and $\mathcal{G}_2$ be two bundle gerbes over $M$. In
a (strictly associative) 2-category, a 1-morphism $\mathcal{A}:\mathcal{G}_1
\to \mathcal{G}_2$ is called invertible or 1-isomorphism, if there is a  1-morphism
$\mathcal{A}^{-1}:\mathcal{G}_2 \to \mathcal{G}_1$  in the opposite direction,
together with 2-isomorphisms $i_{l}:\mathcal{A}^{-1}
\circ \mathcal{A}\Rightarrow \id_{\mathcal{G}_1}$ and $i_r:\id_{\mathcal{G}_2}
\Rightarrow\mathcal{A} \circ \mathcal{A}^{-1}$ such that the diagram
\begin{equation}
\label{44}
\alxydim{@=1.5cm}{\mathcal{A} \circ \mathcal{A}^{-1} \circ \mathcal{A} \ar@{<=}[d]_-{i_r
\circ \id_{\mathcal{A}}} \ar@{=>}[r]^-{\id_{\mathcal{A}}
\circ i_l}  & \mathcal{A} \circ \id_{\mathcal{G}_1} \ar@{=>}[d]^-{\lambda_{\mathcal{A}}}
\\ \id_{\mathcal{G}_2} \circ\mathcal{A} \ar@{=>}[r]_-{\rho_{\mathcal{A}}}
& \mathcal{A}}
\end{equation}
of 2-isomorphisms is commutative. The inverse 1-isomorphism $\mathcal{A}^{-1}$ is unique up to
isomorphism.

Notice that if $\beta:\mathcal{A} \Rightarrow \mathcal{A}'$ is  a 2-morphism
between invertible 1-morphisms we can form a 2-morphism $\beta^{\#}:
\mathcal{A}'^{-1} \Rightarrow \mathcal{A}^{-1}$ using the 2-isomorphisms
$i_r$ for $\mathcal{A}^{-1}$ and $i_l$ for $\mathcal{A}'^{-1}$. Then, diagram
(\ref{44}) induces the equation $\id_{\mathcal{A}}^{\#} = \id_{\mathcal{A}^{-1}}$.

\begin{proposition}
\label{prop2}
A 1-morphism $\mathcal{A}:\mathcal{G}_1 \to \mathcal{G}_2$ in $\mathfrak{BGrb}(M)$ is invertible
if and only if the vector bundle of $\mathcal{A}$ is of rank $1$.
\end{proposition}

\proof
Suppose that $\mathcal{A}$ is invertible, and let $n$ be the rank of its
vector bundle. Let $\mathcal{A}^{-1}$ be an inverse 1-morphism with a vector
bundle of rank $m$. By definition, the composed 1-morphisms $\mathcal{A}\circ
\mathcal{A}^{-1}$ and $\mathcal{A}^{-1} \circ \mathcal{A}$ have vector bundles
of rank $nm$, which has -- to admit the existence of the 2-isomorphisms
$i_l$ and $i_r$ -- to coincide with the rank of the vector bundle of
the identity 1-morphisms $\id_{\mathcal{G}_1}$ and $\id_{\mathcal{G}_2}$ respectively, which is 1. So $n=m=1$. The other inclusion is shown  below
by an explicit construction of an inverse 1-morphism $\mathcal{A}^{-1}$ to a 1-morphism $\mathcal{A}$ with vector bundle of rank 1.
\endofproof

Let a 1-morphism $\mathcal{A}$ consist of a surjective submersion $\zeta:Z
\to Y_1 \times_M Y_2$, of a line bundle $A$ over $Z$ and of an isomorphism $\alpha$ of line
bundles over $Z \times_M Z$. We explicitly
construct an inverse 1-morphism $\mathcal{A}^{-1}$: it has the  surjective
submersion $Z \to Y_1 \times_M Y_2 \to Y_2 \times_M Y_1$, where the first
map
is $\zeta$ and the second one exchanges the factors, the dual line bundle $A^{*}$ over $Z$ and the isomorphism
\begin{equation}
\alxydim{@C=0.5cm}{L_2 \otimes \zeta_2^{*}A^{*} \ar@{=}[r] & \zeta_1^{*}A^{*} \otimes \zeta_1^{*}A \otimes L_2 \otimes
\zeta_2^{*}A^{*} \ar[d]^{\id \otimes \alpha^{-1}
\otimes \id} & \\ & \zeta_1^{*}A^{*}\otimes L_1 \otimes \zeta_2^{*}A \otimes
\zeta_2^{*}A^{*} \ar@{=}[r] & \zeta_1^{*}A^{*}
\otimes L_1.}
\end{equation} 
 Axiom \hyperlink{1M1}{(1M1)} for the 1-morphism $\mathcal{A}^{-1}$ is satisfied because $A^{*}$ has the negative
curvature, and axiom \hyperlink{1M2}{(1M2)}  follows from the one for $\mathcal{A}$.

 To construct
the 2-isomorphism $i_l:\mathcal{A}^{-1}\circ \mathcal{A} \Rightarrow \id_{\mathcal{G}_1}$
notice that the 1-morphism $\mathcal{A}^{-1}\circ \mathcal{A}$ consists of
the  line bundle
$\zeta_1^{*}A\otimes \zeta_2^{*}A^{*}$ over $\tilde Z=Z \times_{Y_2} Z$.
We identify $\tilde Z \cong \tilde Z \times_{P} Y_1^{[2]}$, where $P=Y_1^{[2]}$,
which allows
us to choose
 a triple $(\tilde Z,\id_{\tilde Z},\beta_{\tilde Z})$ defining $i_l$. In
 this triple, the isomorphism
 $\beta_{\tilde Z}$ is defined to be the  composition
\begin{equation}
\label{23}
\hspace{-1cm}
\alxydim{@C=0.4cm}{\zeta_1^{*}A \otimes \zeta_2^{*}A^{*} \ar[rrr]^-{\id \otimes t_{\mu_2}^{-1}
\otimes \id} &&& \zeta_1^{*}A \otimes
\Delta^{*}L_2 \otimes \zeta_2^{*}A^{*} \ar[rr]^-{\alpha^{-1} \otimes \id} &&  L_1 \otimes \zeta_2^{*}A
\otimes \zeta_2^{*}A^{*} = L_1.}
\hspace{-1cm}
\end{equation}
Axiom \hyperlink{2M}{(2M)} for the isomorphism $\beta_{\tilde Z}$ follows from axiom (1M2) of $\mathcal{A}$, so that the triple $(\tilde Z,\id_{\tilde Z},\beta_{\tilde Z})$ defines a
2-isomorphism $i_l: \mathcal{A}^{-1}\circ \mathcal{A} \Rightarrow \id_{\mathcal{G}_1}$.
The 2-isomorphism $i_r: \id_{\mathcal{G}_2} \Rightarrow \mathcal{A} \circ \mathcal{A}^{-1}$
is constructed analogously: here we take the isomorphism
\begin{equation}
\label{24}
\hspace{-1cm}
\alxydim{@C=0.4cm}
{
L_2 =\zeta_1^{*}A^{*} \otimes \zeta_1^{*}A \otimes L_2 \ar[rr]^-{\id \otimes \alpha^{-1}} 
&&  
\zeta_1^{*}A^{*} \otimes
\Delta^{*}L_1 \otimes \zeta_2^{*}A \ar[rrr]^-{\id \otimes t_{\mu_1}
\otimes \id} 
&&& \zeta_1^{*}A^{*} \otimes \zeta_2^{*}A.
}
\hspace{-1cm}
\end{equation}
of line bundles over $W$. Notice that by using the pairing $A^{*}\otimes
A=1$ we have used that $A$ is a line bundle as assumed. Finally, the commutativity
of diagram (\ref{44}) follows from axiom
\hyperlink{1M2}{(1M2)} of $\mathcal{A}$. 

\medskip

Proposition \ref{prop2} shows that we have many 1-morphisms in $\mathfrak{BGrb}(M)$
which are not invertible, in contrast to the 2-groupoid of bundle gerbes defined
in \cite{stevenson1}. Notice that we have already benefited from the simple
definition of the composition $\mathcal{A}^{-1} \circ \mathcal{A}$, which
makes it also easy to see that it is compatible with the construction of
 inverse 1-morphisms $\mathcal{A}^{-1}$:
\begin{equation}
(\mathcal{A}_2 \circ \mathcal{A}_1)^{-1} = \mathcal{A}_1^{-1} \circ \mathcal{A}_2^{-1}
\text{.}
\end{equation}

\medskip

\subsection{Additional Structures}

\label{sec8}

The 2-category of bundle gerbes has  natural definitions of pullbacks,
tensor products and dualities; all of them have been introduced for objects in \cite{murray,murray2}.

Pullbacks and tensor products of 1-morphisms and 2-morphisms can  also be
defined in a natural way, and we do not carry out the details here. Summarizing, 
the monoidal structure on  $\mathfrak{BGrb}(M)$ is a strict 2-functor
\begin{equation}
\otimes \lmaps \mathfrak{BGrb}(M) \times \mathfrak{BGrb}(M) \lto \mathfrak{BGrb}(M)\text{,}
\end{equation}
for which the trivial bundle gerbe $\mathcal{I}_0$ is a strict  tensor unit,
i.e. 
\begin{equation}
\mathcal{I}_0 \otimes \mathcal{G} = \mathcal{G}=\mathcal{G} \otimes
\mathcal{I}_0\text{.}
\end{equation}
The idea of the definition
of $\otimes$ is to take fibre products of the involved surjective submersions,
to pull back all the structure to this fibre product and then to use the
monoidal structure of the category of  vector bundles over that space. This was assumed
to be strict, and so is $\otimes$.
Pullbacks for the 2-category $\mathfrak{BGrb}(M)$ are implemented by    strict 
monoidal  2-functors
\begin{equation}
f^{*}\lmaps\mathfrak{BGrb}(M) \lto \mathfrak{BGrb}(X)
\end{equation}
associated to every smooth map $f:X \to M$ in the way that
$g^{*} \circ f^{*}=(f \circ g)^{*}$
for a second smooth map $g:Y \to X$. The idea of its definition is, to pull
back surjective submersions, for instance
\begin{equation}
\alxy{f^{-1}Y \ar[r]^-{\tilde f} \ar[d] &Y \ar[d]^{\pi} \\ X \ar[r]_{f} &
M}
\end{equation}
and then pull back the structure over $Y$ along the covering map $\tilde
f$. The 2-functors $\otimes$ and $f^{*}$ are all compatible with
the assignment of inverses $\mathcal{A}^{-1}$ to 1-morphisms $\mathcal{A}$
from subsection \ref{sec6}:
\begin{equation}
f^{*}(\mathcal{A}^{-1})=(f^{*}\mathcal{A})^{-1}
\quad\text{ and }\quad
(\mathcal{A}_1 \otimes \mathcal{A}_2)^{-1} = \mathcal{A}_1^{-1} \otimes \mathcal{A}_2^{-1}\text{.}
\end{equation}
Also the trivial bundle gerbes $\mathcal{I}_{\rho}$ behave naturally under pullbacks
and tensor products:
\begin{equation}
f^{*}\mathcal{I}_{\rho} = \mathcal{I}_{f^{*}\rho}
\quad\text{ and }\quad
\mathcal{I}_{\rho_1} \otimes \mathcal{I}_{\rho_2} = \mathcal{I}_{\rho_1 +
\rho_2}\text{.}
\end{equation}

\medskip

To define a duality we are a bit more precise, because this has yet not been
done systematically in the literature. Even though we will strictly concentrate
on what we  need in section \ref{sec7}. For those purposes,
it is enough to understand the duality
as a strict 2-functor
\begin{equation}
()^{*}: \mathfrak{BGrb}(M)^{\mathrm{op}} \to \mathfrak{BGrb}(M)
\end{equation}
where the opposed 2-category $\mathfrak{BGrb}(M)^{\mathrm{op}}$
has all 1-morphisms reversed, while
the 2-morphisms are as before. This
2-functor will satisfy the identity
\begin{equation}
()^{**}=\id_{\mathfrak{BGrb}(M)}\text{.}
\end{equation}

We now give the complete definition
of the functor $()^{*}$ on objects,
1-morphisms and 2-morphisms. For a given bundle gerbe $\mathcal{G}$,
the dual bundle gerbe $\mathcal{G}^{*}$ consists of the same surjective submersion
$\pi:Y \to M$, the 2-form $-C \in \Omega^2(Y)$, the line bundle $L^{*}$
over $Y^{[2]}$ and the  isomorphism 
\begin{equation}
\mu^{*-1}:\pi_{12}^{*}L^{*}\otimes \pi_{23}^{*}L^{*}\to
\pi_{13}^{*}L^{*}
\end{equation}
of line bundles over $Y^{[3]}$. This structure clearly satisfies the axioms of a bundle
gerbe. We obtain immediately
\begin{equation}
\label{81}
\mathcal{G}^{**}=\mathcal{G}
\quad\text{ and }\quad
(\mathcal{G} \otimes \mathcal{H})^{*}
= \mathcal{H}^{*} \otimes \mathcal{G}^{*}
\text{,}
\end{equation}
and for the trivial bundle gerbe $\mathcal{I}_{\rho}$  the equation
\begin{equation}
\mathcal{I}_{\rho}^{*}=\mathcal{I}_{-\rho}\text{.}
\end{equation}

For a 1-morphisms $\mathcal{A}:\mathcal{G}_1
\to \mathcal{G}_2$ consisting of
a vector bundle $A$ over $Z$ with
surjective submersion $\zeta: Z
\to P$ with $P:=Y_1 \times_M Y_2$ and of an
isomorphism $\alpha$ of vector
bundles over $Z \times_M
Z$, we define the dual 1-morphism
\begin{equation}
\mathcal{A}^{*}\lmaps\mathcal{G}_2^{*}
\lto \mathcal{G}_1^{*}
\end{equation}
as follows: its surjective submersion is the pullback of $\zeta$
along the exchange map $s : P' \to P$, with $P':=Y_2 \times Y_1$; that is
a surjective submersion $\zeta':Z' \to P'$ and a covering map $s_Z$ in the commutative
diagram
\begin{equation}
\label{83}
\alxydim{@C1.2cm@R=1.2cm}{Z' \ar[d]_{\zeta'} \ar[r]^-{s_Z} & Z \ar[d]^{\zeta} \\ P' \ar[r]_{s} & P\text{.}}
\end{equation}
The  vector
bundle of $\mathcal{A}^{*}$ is  $A':=s_Z^{*}A$ over $Z'$ and  its isomorphism
is
\begin{equation}
\label{90}
\alxy{L_2^{*} \otimes \zeta_2^{\prime*}A'
\ar@{=}[r] & L_2^{*} \otimes L_1
\otimes \zeta_2^{\prime*}s_Z^{*}A
\otimes L_1^{*} \ar[d]^{\id \otimes \tilde s^{*}
\alpha \otimes \id} & \\ & L_2^{*}
\otimes \zeta_1^{\prime*}s_Z^{*}A \otimes L_2
\otimes L_1^{*} \ar@{=}[r]  & \zeta_1^{\prime*}A'
\otimes L_1^{*}\text{.}}
\end{equation} 
Axiom \hyperlink{1M1}{(1M1)}
is satisfied since the dual bundle
gerbes have 2-forms with opposite
signs, 
\begin{equation}
\mathrm{curv}(A')= s_Z^{*}\mathrm{curv}(A) =s_Z^{*}(C_2 - C_1 )=C_2 - C_1 = (-C_1)
- (-C_2)\text{.}
\end{equation}   
Axiom \hyperlink{1M2}{(1M2)}
relates the isomorphism (\ref{90})
to the isomorphisms $\mu_1^{*-1}$
and $\mu_2^{*-1}$ of the dual bundle
gerbes. It can be deduced from
axiom \hyperlink{1M2}{(1M2)} of
$\mathcal{A}$ using the following
general fact, applied to $\mu_1^{*}$
and $\mu_2^{*}$: 
the dual $f^{*}$ of an isomorphism $f:L_1 \to L_2$
of line bundles coincides
with the isomorphism
\begin{equation}
\label{75}
\alxydim{@C=2cm}{L_2^{*} =L_2^{*} \otimes L_1 \otimes L_1^{*} \ar[r]^-{\id \otimes f
\otimes \id} & L_2^{*} \otimes L_2 \otimes L_1^{*} = L_1^{*}\text{,}}
\end{equation}
defined using the duality on line bundles. 

Dual 1-morphisms defined like
this have the properties
\begin{equation}
\label{82}
\mathcal{A}^{**} = \mathcal{A}
\quad\text{, }\quad
(\mathcal{A}' \circ \mathcal{A})^{*}
= \mathcal{A}^{*}\circ \mathcal{A}'^{*}
\quad\text{ and }\quad
(\mathcal{A}_1 \otimes \mathcal{A}_2)^{*} = \mathcal{A}_2^{*} \otimes \mathcal{A}_1^{*}\text{.}
\end{equation} 

Finally, for a 2-morphism $\beta: \mathcal{A}_1 \Rightarrow \mathcal{A}_2$ we
define the dual 2-morphism
\begin{equation}
\beta^{*}\lmaps \mathcal{A}_1^{*} \Longrightarrow \mathcal{A}_2^{*}
\end{equation}
in the following way. If $\beta$ is represented by a triple $(W,\omega,\beta_W)$
with an isomorphism $\beta_W: A_1 \to A_2$ of vector bundles over $W$, we consider
the pullback of $\omega: W \to Z_1 \times_P Z_2$ along $s_{Z_1} \times s_{Z_2}:
Z_1' \times_{P'} Z_2' \to Z_1 \times_P Z_2$, where $Z_1$, $Z_2'$ and $P'$ are
as in (\ref{83}), and $s_{Z_1}$ and $s_{Z_2}$ are the respective covering maps. This
gives a commutative diagram
\begin{equation}
\alxydim{@C2cm@R=1.2cm}{W' \ar[d]_{\omega'} \ar[r]^-{s_W} & W \ar[d]^{\omega} \\ Z_1' \times_{P'} Z_2' \ar[r]_{s_{Z_1} \times s_{Z_2}} & Z_1 \times_P Z_2\text{.}}
\end{equation}
Now consider the triple $(W',\omega',s_W^{*} \beta_W)$
with the isomorphism 
\begin{equation}
s_W^{*}\beta_W\lmaps s_{Z_1}^{*}A_1 \lto s_{Z_2}^{*}A_2
\end{equation}
of vector bundles over $W'$. It satisfies axiom \hyperlink{2M}{(2M)}, and
its equivalence class does not depend on the choice of the representative
of $\beta$. So we define the dual 2-morphism $\beta^{*}$ to be this class.  Dual 2-morphisms are
compatible with vertical and horizontal compositions
\begin{equation}
\label{85}
(\beta_2 \circ \beta_1)^{*} = \beta_1^{*} \circ \beta_2^{*}
\quad\text{ and }\quad
(\beta \bullet \beta')^{*} = \beta^{*} \bullet \beta'^{*}
\end{equation}
and satisfy furthermore
\begin{equation}
\label{32}
\beta^{**}=\beta
\quad\text{ and }\quad
(\beta_1 \otimes \beta_2)^{*} = \beta_2^{*} \otimes \beta_1^{*} \text{.}
\end{equation}

We can use adjoint 2-morphisms
in the following situation: if $\mathcal{A}:\mathcal{G}
\to \mathcal{H}$ is an invertible
1-morphism with inverse $\mathcal{A}^{-1}$ and associated  2-isomorphisms $i_l: \mathcal{A}^{-1}
\circ \mathcal{A} \Rightarrow \id_{\mathcal{G}}$ and $i_r:\id_{\mathcal{H}}
\Rightarrow \mathcal{A} \circ \mathcal{A}
^{-1}$, their duals  $i^{*}_l$ and $i^{*}_r$  show that $(\mathcal{A}^{-1})^{*}$ is an inverse
to $\mathcal{A}^{*}$. Since
inverses are unique up to isomorphism,
\begin{equation}
(\mathcal{A}^{*})^{-1} \cong (\mathcal{A}^{-1})^{*}\text{.}
\end{equation}

Summarizing, equations (\ref{81}), (\ref{82}), (\ref{85}) and (\ref{32}) show that $()^{*}$
is a monoidal strict 2-functor, which is strictly involutive. Let us finally
mention that it is also compatible with pullbacks:
\begin{equation}
f^{*}(\mathcal{G}^{*}) = (f^{*}\mathcal{G})^{*}
\quad\text{, }\quad
f^{*}\mathcal{A}^{*} = (f^{*}\mathcal{A})^{*}
\quad\text{ and }\quad
f^{*}\beta^{*} = (f^{*}\beta)^{*}\text{.}
\end{equation}

\section{Descent Theory for Morphisms}

\label{sec5}

In this section we compare 1-morphisms between bundle gerbes in the sense
of  Definition \ref{def1} with 1-morphisms whose surjective submersion
$\zeta:Z \to Y_1 \times_M Y_2$ is the identity, like in \cite{stevenson1}. For this purpose, we introduce the subcategory  $\mathfrak{Hom}_{FP}(\mathcal{G}_1,\mathcal{G}_2)$
of the morphism category $\mathfrak{Hom}(\mathcal{G}_1,\mathcal{G}_2)$, where all smooth manifolds $Z$ and $W$ appearing in the definitions of 1- and 2-morphisms
are equal to the fibre product $P:=Y_1 \times_M Y_2$.
Explicitly, an object in $\mathfrak{Hom}_{FP}(\mathcal{G}_1,\mathcal{G}_2)$
is a 1-morphism $\mathcal{A}:\mathcal{G}_1 \to \mathcal{G}_2$
whose surjective submersion is the identity $\id_P$ and a morphism in $\mathfrak{Hom}_{FP}(\mathcal{G}_1,\mathcal{G}_2)$
is a 2-morphism $\beta:\mathcal{A}_1 \Rightarrow
\mathcal{A}_2$ which can be represented by the triple $(P,\omega,\beta)$
where $\omega: P \to P \times_P P \cong P$ is the identity. 
\begin{theorem} 
\label{th1}
The inclusion functor
\begin{equation*}
D\lmaps\mathfrak{Hom}_{FP}(\mathcal{G}_1,\mathcal{G}_2)
\lto \mathfrak{Hom}(\mathcal{G}_1,\mathcal{G}_2)
\end{equation*}
is an equivalence of categories.
\end{theorem}

In the proof we will
make use of the fact that vector bundles form a stack, i.e.
 fibred category satisfying a gluing condition. To
make this gluing condition concrete, we define for a surjective submersion $\zeta: Z \to P$ a category $\mathfrak{Des}(\zeta)$
as follows. Its objects  are pairs $(A,d)$, where
$A$ is a vector bundle over $Z$ and 
\begin{equation}
d\lmaps\zeta_1^{*}A \lto \zeta_2^{*}A
\end{equation}
is an isomorphism of vector bundles over $Z^{[2]}$ 
such that
\begin{equation}
\label{30}
\zeta_{13}^{*}d = \zeta_{23}^{*}d \circ \zeta_{12}^{*}d\text{.}
\end{equation}
A morphism $\alpha: (A,d) \to (A', d')$ in $\mathfrak{Des}(\zeta)$ is an isomorphism $\alpha: A \to
A'$ of vector bundles over $Z$ such
that the diagram
\begin{equation}
\alxydim{@C=1.4cm}{\zeta_1^{*}A \ar[d]_{d} \ar[r]^{\zeta_1^{*}\alpha} & \zeta_1^{*}A'  \ar[d]^{d'} \\ \zeta_2^{*}A \ar[r]_{\zeta_2^{*}\alpha} & \zeta_2^{*}A'}
\label{descmorph}\end{equation}
of isomorphisms of vector bundles over $Z^{[2]}$ is commutative.
The composition of morphisms is just the
composition of isomorphisms of vector bundles.
Now, the gluing condition for the stack of vector bundles is that the pullback along $\zeta$ is an equivalence
\begin{equation}
\zeta^{*}\lmaps \mathfrak{Bun}(P) \lto \mathfrak{Des}(\zeta)
\end{equation}
between the category $\mathfrak{Bun}(P)$ of vector bundles
over $P$ and the category $\mathfrak{Des}(\zeta)$. 

\medskip

Proof of Theorem \ref{th1}. We show that the faithful functor $D$ is an equivalence of categories by proving 
 (a) that it is essentially surjective and (b) that the subcategory $\mathfrak{Hom}_{FP}(\mathcal{G}_1,\mathcal{G}_2)$
 is full. 

For (a) we have to show that for every 1-morphism $\mathcal{A}:\mathcal{G}_1
\to \mathcal{G}_2$ with arbitrary surjective submersion $\zeta:Z \to P$ there is an isomorphic 1-morphism $\mathcal{S}_{\mathcal{A}}:\mathcal{G}_1 \to \mathcal{G}_2$
with surjective submersion $\id_P$.
Notice that the isomorphism $\mathrm{d}_\mathcal{A}:\zeta_1^{*}A \to
\zeta_2^{*}A$ of vector bundles over $Z^{[2]}$ from Lemma \ref{lem3} satisfies the cocycle condition (\ref{30}), so that $(A,\mathrm{d}_{\mathcal{A}})$ is an object in $\mathfrak{Des}(\zeta)$.
Now consider the surjective submersion $\zeta^2: Z \times_M Z\to P^{[2]}$.
By Lemma \ref{lem3} b) and under the identification of $Z^{[2]} \times_M
Z^{[2]}$ with $(Z \times_M Z) \times_{P^{[2]}} (Z \times_M Z)$ the
 diagram
 \begin{equation}
\alxydim{@C=1.5cm}{L_1 \otimes \zeta_2^{*}A \ar[d]_{1 \otimes \zeta_{24}^{*}\mathrm{d}_{\mathcal{A}}} \ar[r]^{\zeta_{12}^{*}\alpha}
& \zeta_1^{*}A \otimes L_2 \ar[d]^{\zeta_{13}^{*}\mathrm{d}_{\mathcal{A}} \otimes 1} \\ L_1 \otimes \zeta_4^{*}A \ar[r]_{\zeta_{34}^{*}\alpha}
& \zeta_3^{*}A \otimes L_2}
\end{equation}   
of isomorphisms of vector bundles over $(Z \times_M Z) \times_{P^{[2]}} (Z \times_M Z)$ is commutative, and shows that $\alpha$ is a morphism in $\mathfrak{Des}(\zeta^2)$.
Now we use that $\zeta^{*}$ is an equivalence of categories: we choose a vector bundle $S$ over $P$ together with an isomorphism $\beta:\zeta^{*}S\to
A$ of vector bundles over $Z$, and an isomorphism
\begin{equation}
\sigma\lmaps L_1 \otimes \zeta_2^{*}S \lto \zeta_1^{*}S \otimes L_2
\end{equation}
 of vector bundles over $P \times_M P$ such that the diagram
\begin{equation}
\label{100}
\alxydim{@C=1.5cm}{L_1 \otimes \zeta_2^{*}\zeta^{*}S \ar[d]_{\id \otimes \zeta_2^{*}\beta} \ar[r]^{\zeta^{*}\sigma} & \zeta_1^{*}\zeta^{*}S \otimes L_2 \ar[d]^{\zeta_1^{*}\beta
\otimes \id}
\\ L_1 \otimes \zeta_2^{*}A  \ar[r]_{\alpha} & \zeta_1^{*}A \otimes L_2}
\end{equation} 
 of isomorphisms of vector bundles over $Z \times_M Z$  is commutative. Since
 $\zeta$ is an equivalence of categories, the axioms of $\mathcal{A}$ imply
 the ones  of the 1-morphism $\mathcal{S}_{\mathcal{A}}$ defined by the surjective submersion
 $\id_P$, the vector bundle $S$ over $P$ and the isomorphism $\sigma$ over
 $P^{[2]}$. Finally, the triple $(Z \times_P P,\id_Z,\beta)$ with $Z \cong Z \times_P
 P$ defines a 2-morphism $\mathcal{S}_{\mathcal{A}}
 \Rightarrow \mathcal{A}$, whose
 axiom \hyperlink{2M}{(2M)} is
 (\ref{100}).

(b) We have to show that  any morphism $\beta:\mathcal{A}\Rightarrow \mathcal{A}'$ in $\mathfrak{Hom}(\mathcal{G}_1,\mathcal{G}_2)$
 between objects $\mathcal{A}$ and $\mathcal{A}'$ in $\mathfrak{Hom}_{FP}(\mathcal{G}_1,\mathcal{G}_2)$ 
is already a morphism in $\mathfrak{Hom}_{FP}(\mathcal{G}_1,\mathcal{G}_2)$.
Let $(W,\omega,\beta_W)$ be any representative of $\beta$ with a surjective
submersion $\omega:W \to P$ and an isomorphism
$\beta_W: \omega^{*}A \to \omega^{*}A'$ of vector bundles over $W$. The restriction
of axiom \hyperlink{2M}{(2M)} for the triple $(W,\omega,\beta_W)$ from $W \times_M W$ to $W
\times_P W$ shows $\omega_1^{*}\beta_W = \omega_2^{*}\beta_W$.  This shows
that $\beta_W$ is a morphism in the descent category $\mathfrak{Des}(\omega)$.
Let $\beta_P: A \to A'$ be an isomorphism of vector bundles over $P$ such
that 
\begin{equation}
\label{31}
\omega^{*}\beta_P=\beta_W
\end{equation}
Because $\omega$
is an equivalence of categories, the triple $(P,\id_P, \beta_P)$ defines
a 2-morphism from $\mathcal{A}$ to $\mathcal{A}'$ being a morphism in $\mathfrak{Hom}_{FP}(\mathcal{G}_1,\mathcal{G}_2)$. Equation (\ref{31}) shows
that the triples $(P,\id_P,\beta_P)$ and $(W,\omega,\beta_W)$ are equivalent.
\endofproof

In the remainder of this section we present two corollaries of Theorem \ref{th1}.
First, and most importantly, we make contact to the notion of a stable isomorphism
between bundle gerbes. By definition \cite{murray2}, a
stable isomorphism is a 1-morphism, whose surjective submersion is the identity
 $\id_P$ on the fibre product of the surjective submersions of the two bundle
 gerbes, and whose vector bundle over $P$ is a line bundle. From
Proposition \ref{prop2} and
Theorem \ref{th1} we obtain  

\begin{corollary}
\label{co1}
There exists a stable isomorphism $\mathcal{A}:\mathcal{G}_1 \to\mathcal{G}_2$
if and only if the bundle gerbes are isomorphic objects in $\mathfrak{BGrb}(M)$. \end{corollary}

It is shown in \cite{murray2}
that the set of stable isomorphism classes of bundle gerbes over $M$ is a
group (in virtue of the monoidal structure) which is isomorphic to the Deligne cohomology
group $\mathrm{H}^2(M,\mathcal{D}(2))$. This is a very important fact which
connects the theory of bundle gerbes to other theories of gerbes, for instance,
to Dixmier-Douady sheaves
of groupoids \cite{brylinski1}. Corollary \ref{co1} states that although
our definition of morphisms differs from the one of \cite{murray2}, the bijection
between isomorphism classes of bundle gerbes and 
the Deligne cohomology group is persistent.

\medskip

Second, Theorem \ref{th1} admits to use existing classification results
for 1-isomorphisms. Consider the full subgroupoid $\mathfrak{Aut}(\mathcal{G})$ of
$\mathfrak{Hom}(\mathcal{G},\mathcal{G})$ associated
to a bundle gerbe $\mathcal{G}$, which consists
of all 1-isomorphisms $\mathcal{A}:\mathcal{G} \to \mathcal{G}$, and all
(necessarily invertible) 2-morphisms between those. From Theorem \ref{th1}
and Lemma 2 of \cite{schreiber1} we obtain

\begin{corollary}
The skeleton of the groupoid $\mathfrak{Aut}(\mathcal{G})$, i.e. the set
of isomorphism classes of 1-isomorphisms $\mathcal{A}:\mathcal{G} \to \mathcal{G}$
is a torsor over the group $\mathrm{Pic}_0(M)$ of isomorphism classes of
flat line bundles over $M$.
\end{corollary} 

In  2-dimensional conformal field
theory, where a bundle gerbe $\mathcal{G}$ is considered to be a part of the background
field, the groupoid $\mathfrak{Aut}(\mathcal{G})$ may be called the groupoid
of gauge transformations of $\mathcal{G}$. The above corollary classifies
such gauge transformation up to equivalence.

\section{Some important Examples of Morphisms}

\label{sec4}

To discuss holonomies of bundle gerbes, it is essential to establish an  equivalence
between the morphism categories between trivial bundle gerbes  over $M$ and vector
bundles of certain curvature over $M$. Given two 2-forms $\rho_1$ and $\rho_2$ on $M$, consider  the category $\mathfrak{Hom}_{FP}(\mathcal{I}_{\rho_1},\mathcal{I}_{\rho_2})$.
An object $\mathcal{A}:\mathcal{I}_{\rho_1} \to \mathcal{I}_{\rho_2}$ consists
of
the smooth manifold $Z = M$ with
the surjective submersion $\zeta=\id_M$, a vector bundle $A$ over $M$ and an isomorphism
$\alpha:A \to A$. Axiom \hyperlink{1M2}{(1M2)} states 
\begin{equation}
\label{200}
\frac{1}{n}\mathrm{tr}(\mathrm{curv}(A)) = \rho_2-\rho_1
\end{equation}
 with $n$ the rank of $A$,
and axiom \hyperlink{1M2}{(1M2)} reduces to $\alpha^2=\alpha$, which in turn means $\alpha=\id_A$.
Together with Theorem \ref{th1}, this defines a canonical equivalence of
categories
\begin{equation}
\mathrm{Bun}\lmaps\mathfrak{Hom}(\mathcal{I}_{\rho_1},\mathcal{I}_{\rho_2}) \lto \mathfrak{Bun}_{\rho_2-\rho_1}(M)\text{,}
\end{equation}
where $\mathfrak{Bun}_{\rho}(M)$ is the category of vector bundles $A$ over $M$
whose curvature satisfies (\ref{200}).
Its following properties can directly be deduced
from the definitions.

\begin{proposition}
\label{prop1}
The functor $\mathrm{Bun}$ respects the structure of the 2-category
of bundle gerbes, namely:    
\begin{itemize}
\item[a)]
the composition of 1-morphisms,
\begin{equation*}
\mathrm{Bun}(\mathcal{A}_2 \circ \mathcal{A}_1)=\mathrm{Bun}(\mathcal{A}_1) \otimes \mathrm{Bun}(\mathcal{A}_2)\quad\text{ and }\quad\mathrm{Bun}(\id_{\mathcal{\mathcal{I}_{\rho}}})=1\text{.}
\end{equation*}

\item[b)]
the assignment of inverses to invertible 1-morphisms,
\begin{equation*}
\mathrm{Bun}(\mathcal{A}^{-1}) = \mathrm{Bun}(\mathcal{A})^{*}\text{.}
\end{equation*}

\item[c)]
the monoidal structure,
\begin{equation*}
\mathrm{Bun}(\mathcal{A}_1 \otimes \mathcal{A}_2) = \mathrm{Bun}(\mathcal{A}_1)
\otimes \mathrm{Bun}(\mathcal{A}_2)\text{.}
\end{equation*}

\item[d)]
pullbacks,
\begin{equation*}
\mathrm{Bun}(f^{*}\mathcal{A})=f^{*}\mathrm{Bun}(\mathcal{A})
\quad\text{ and }\quad
\mathrm{Bun}(f^{*}\beta)=f^{*}\mathrm{Bun}(\beta)\text{.}
\end{equation*}

\item[e)] the duality\begin{equation*}
\mathrm{Bun}(\mathcal{A}^{*}) = \mathrm{Bun}(\mathcal{A})
\quad\text{ and }\quad 
\mathrm{Bun}(\beta^{*}) = \mathrm{Bun}(\beta)\text{.}
\end{equation*}

\end{itemize}
\end{proposition}

In the following subsections we see how the 2-categorial structure of bundle gerbes and the functor $\mathrm{Bun}$ can
be used to give natural  definitions of surface holonomy in several
situations.

\subsection{Trivializations}

We give the following natural definition of a trivialization.

\begin{definition}
A trivialization of a bundle gerbe $\mathcal{G}$ is a 1-isomorphism
\begin{equation*}
\mathcal{T}\lmaps \mathcal{G} \lto \mathcal{I}_{\rho}\text{.}
\end{equation*}
\end{definition}

Let us write out the details of
such a 1-isomorphism. By Theorem
\ref{th1} we may assume that the
surjective submersion of $\mathcal{T}$ is the identity $\id_P$ on $P:=Y
\times_M M\cong Y$ with projection $\pi$ to $M$.
Then, $\mathcal{T}$ consists further
of a line bundle $T$ over $Y$,
and of an isomorphism $\tau: L
\otimes \pi_2^{*}T \to \pi_1^{*}T$
of line bundles over $Y^{[2]}$.
Axiom \hyperlink{1M2}{(1M2)} gives $\pi_{13}^{*}\tau
\circ (\mu \otimes \id) = \pi_{12}^{*}\tau
\circ \pi_{23}^{*}\tau$ as an equation
of isomorphisms of line bundles
over $Y^{[3]}$. This is exactly
the definition of a trivialization
one finds in the literature \cite{carey2}.
Additionally, axiom \hyperlink{1M2}{(1M2)} gives
$\mathrm{curv}(T) = \pi^{*}\rho - C$. If one specifies $\rho$ not as a part
of the definition of a trivialization, it is uniquely determined by this
equation. 

\medskip

\medskip

Trivializations are essential for the definition of holonomy around closed
oriented surfaces.

\begin{definition}
\label{def6}
Let $\phi:\Sigma \to M$ be a smooth map from a closed oriented surface $\Sigma$ to a smooth manifold $M$,
and let $\mathcal{G}$ a bundle gerbe
over $M$. Let
\begin{equation*}
\mathcal{T}\lmaps \phi^{*}\mathcal{G} \lto \mathcal{I}_{\rho}
\end{equation*}
be any trivialization. The holonomy of $\mathcal{G}$ around $\phi$ is defined
as
\begin{equation*}
\mathrm{hol}_{\mathcal{G}}(\phi) := \exp \left ( \mathrm{i} \int_{\Sigma} \rho
\right ) \in U(1)\text{.}
\end{equation*}
\end{definition}

In this situation,  the functor $\mathrm{Bun}$ is a powerful tool to prove that this definition does  not depend on the choice of the trivialization:
if $\mathcal{T}':\phi^{*}\mathcal{G} \to \mathcal{I}_{\rho'}$ is another trivialization, the composition
$\mathcal{T} \circ \mathcal{T}^{\prime-1}:\mathcal{I}_{\rho'} \to \mathcal{I}_{\rho}$
corresponds to a line bundle $\mathrm{Bun}(\mathcal{T} \circ \mathcal{T}^{\prime-1})$ over $M$ with curvature $\rho - \rho'$. In particular,
the difference between any two 2-forms $\rho$ is a closed 2-form with integer
periods and vanishes under the exponentiation in the definition
of $\mathrm{hol}_{\mathcal{G}}(\phi)$.

\subsection{Bundle Gerbe Modules}

For oriented surfaces with boundary one has to choose additional structure
on the boundary to obtain a well-defined holonomy. This additional structure
is provided by a vector bundle twisted by the bundle gerbe $\mathcal{G}$
\cite{gawedzki4}, also known as
a $\mathcal{G}$-module.
In our formulation, its definition takes the following form:
\begin{definition}
\label{def4}
Let $\mathcal{G}$ be a bundle gerbe over $M$. A left $\mathcal{G}$-module
is a 1-morphism $\mathcal{E}:
\mathcal{G} \to \mathcal{I}_{\omega}$, and a right $\mathcal{G}$-module is
a 1-morphism $\mathcal{F}:\mathcal{I}_{\omega} \to \mathcal{G}$\text{.}
\end{definition}

Let us compare this definition with the original definition
of (left) bundle gerbe modules in $\cite{bouwknegt1}$. Assume
-- again by Theorem \ref{th1} -- that a left
$\mathcal{G}$-module $\mathcal{E}:\mathcal{G} \to \mathcal{I}_{\omega}$ has
the surjective submersion $\id_P$ with $P\cong Y$. 
Then, it consists of a vector bundle $E$ over $Y$ and of an isomorphism $\epsilon:
L \otimes \pi_2^{*}E \to \pi_1^{*}E$ of vector bundles over $Y^{[2]}$ which satisfies
\begin{equation}
\pi_{13}^{*}\epsilon\circ (\mu \otimes \id)= \pi_{23}^{*}\epsilon \circ \pi_{12}^{*}\epsilon
\end{equation}
by axiom \hyperlink{1M2}{(1M2)}. The curvature of $E$ is restricted by axiom \hyperlink{1M2}{(1M2)} to 
\begin{equation}
\frac{1}{n}\mathrm{tr}(\mathrm{curv}(E)) = \pi^{*}\omega - C
\end{equation}
 with $n$ the rank of $E$. 

\medskip

The definition of bundle gerbe modules as 1-morphisms makes clear that left
and right $\mathcal{G}$-modules form categories $\mathfrak{LMod}(\mathcal{G})$
and $\mathfrak{RMod}(\mathcal{G})$. This is useful for
example to see 
that a 1-isomorphism  $\mathcal{A}:\mathcal{G} \to \mathcal{G}'$ defines
equivalences of categories 
\begin{equation}
\label{99}
\mathfrak{LMod}(\mathcal{G}) \cong \mathfrak{LMod}(\mathcal{G}')
\quad\text{ and }\quad
\mathfrak{RMod}(\mathcal{G}) \cong \mathfrak{RMod}(\mathcal{G}')
\end{equation}
and that there are equivalences between
left modules of $\mathcal{G}$ 
and right modules of $\mathcal{G}^{*}$
(and vice versa), by taking  duals of the respective  1-morphisms. Moreover,
for a trivial bundle gerbe $\mathcal{I}_{\rho}$ the categories $\mathfrak{LMod}(\mathcal{I}_{\rho})$
and $\mathfrak{RMod}(\mathcal{I}_{\rho})$ become canonically equivalent to the category
$\mathfrak{Bun}(M)$ of vector bundles over $M$ via the functor $\mathrm{Bun}$.
We can combine this result with the equivalences (\ref{99})
applied to a trivialization $\mathcal{T}: \mathcal{G} \to \mathcal{I}_{\rho}$
of a bundle gerbe $\mathcal{G}$ over $M$. In detail, a left $\mathcal{G}$-module $\mathcal{E}:\mathcal{G} \to \mathcal{I}_{\omega}$ first becomes a left $\mathcal{I}_{\rho}$-module
\begin{equation}
\mathcal{E}\circ
\mathcal{T}^{-1}\lmaps \mathcal{I}_{\rho} \lto \mathcal{I}_{\omega}
\end{equation}
 which in turn defines the vector bundle $E:=\mathrm{Bun}(\mathcal{E}\circ\mathcal{T}^{-1})$ over $M$. The same applies to right $\mathcal{G}$-modules
$\mathcal{F}:\mathcal{I}_{\omega}\to \mathcal{G}$ which defines a vector bundle $\bar E := \mathrm{Bun}(\mathcal{T}
\circ \mathcal{F})$ over $M$. 

\medskip

A D-brane for the bundle gerbe $\mathcal{G}$ is a submanifold $Q$ of $M$
together with a left $\mathcal{G}|_{Q}$-module. Here $\mathcal{G}|_{Q}$
means the pullback of $\mathcal{G}$ along the inclusion $Q \hookrightarrow
M$. 

\begin{definition}
Let $\mathcal{G}$ be a bundle gerbe over $M$ with D-brane $(Q,\mathcal{E})$
and let $\phi:\Sigma \to M$ be a smooth map from a compact oriented surface
$\Sigma$ with boundary to $M$, such that $\phi(\partial\Sigma) \subset Q$. Let 
\begin{equation*}
\mathcal{T}\lmaps \phi^{*}\mathcal{G} \lto \mathcal{I}_{\rho}
\end{equation*}
be any trivialization of the pullback bundle gerbe $\phi^{*}\mathcal{G}$ and let 
\begin{equation}
E := \mathrm{Bun}(\phi^{*}\mathcal{E}
\circ \mathcal{T}^{-1})
\end{equation}
be the associated vector bundle over
$\partial \Sigma$. The holonomy of $\mathcal{G}$
around $\phi$ is defined as
\begin{equation*}
\mathrm{hol}_{\mathcal{G},\mathcal{E}}(\phi) := \exp \left ( \mathrm{i} \int_{\Sigma}
\rho \right ) \cdot \mathrm{tr}\left ( \mathrm{hol}_E(\partial \Sigma) \right
) \in \C\text{.}
\end{equation*}
\end{definition}

\medskip

The definition does not depend on the choice of the trivialization: for another
trivialization $\mathcal{T}':\phi^{*}\mathcal{G} \to \mathcal{I}_{\rho'}$
and the respective vector bundle $E':=\mathrm{Bun}(\mathcal{E} \circ \mathcal{T}'^{-1})$
we find by Proposition \ref{prop1} a)
\begin{equation}
E' = \mathrm{Bun}(\mathcal{E} \circ \mathcal{T}^{\prime-1}) \cong \mathrm{Bun}(\mathcal{E}
\circ \mathcal{T}^{-1} \circ \mathcal{T} \circ \mathcal{T}^{\prime-1}) =
E \otimes \mathrm{Bun}(\mathcal{T} \circ \mathcal{T}^{\prime-1}).
\end{equation}
Because isomorphic vector bundles have the same holonomies, and the line bundle $\mathrm{Bun}(\mathcal{T} \circ \mathcal{T}'^{-1})$
has curvature $\rho-\rho'$ we obtain
\begin{equation}
\mathrm{tr}\left ( \mathrm{hol}_{E'}(\partial \Sigma) \right
) = \mathrm{tr}\left ( \mathrm{hol}_E(\partial \Sigma) \right
) \cdot \exp \left ( \mathrm{i} \int _{\Sigma} \rho - \rho' \right )\text{.}
\end{equation}
This shows the independence of the choice of the trivialization.

\subsection{Jandl Structures}

\label{sec7}

In this last section, we use the  duality  on the 2-category
$\mathfrak{BGrb}(M)$ introduced in section \ref{sec8} to define the 
holonomy of a bundle gerbe around unoriented, and even unorientable surfaces (without boundary).
For this purpose,  we explain
the concept of a Jandl structure on a bundle gerbe $\mathcal{G}$, which has
been introduced in \cite{schreiber1}, in terms of 1- and 2-isomorphisms of
the 2-category $\mathfrak{BGrb}(M)$.
\begin{definition}
\label{def5}
A Jandl structure $\mathcal{J}$  on a bundle gerbe
$\mathcal{G}$ over $M$ is a collection $(k,\mathcal{A},\varphi)$ of an involution
$k:M \to M$, i.e. a diffeomorphism with $k \circ k=\id_M$, a 1-isomorphism
\begin{equation*}
\mathcal{A}\lmaps k^{*}\mathcal{G} \lto \mathcal{G}^{*}
\end{equation*}
and a 2-isomorphism
\begin{equation*}
\varphi\lmaps k^{*}\mathcal{A} \Longrightarrow \mathcal{A}^{*}
\end{equation*}
which satisfies the condition
\begin{equation*}
k^{*}\varphi = \varphi^{*-1}\text{.}
\end{equation*}
\end{definition}

Notice that the existence of the 2-isomorphism $\varphi$ is only possible because $\mathcal{G}^{**}=\mathcal{G}$
from (\ref{81}), and that the equation $k^{*}\varphi = \varphi^{*-1}$ only makes
sense because $\mathcal{A}^{**}=\mathcal{A}$ from (\ref{82}). Let us now
discuss the relation between Definition \ref{def5} and the original definition of a Jandl structure from \cite{schreiber1}. For this purpose we elaborate
the details. We denote the pullback of
the surjective submersion $\pi:Y \to M$ along $k$ by $\pi_k:Y_k \to M$; for
simplicity we take $Y_k:=Y$ and $\pi_k := k\circ \pi$. Now, we assume by
 Theorem \ref{th1} that the 1-isomorphism $\mathcal{A}$
consists of a line bundle $A$ over $Y_k \times_M Y$. As smooth manifolds,
we can identify $Y_k \times_M Y$ with $P:=Y^{[2]}$; to have an identification
as smooth manifolds with  surjective submersions to $M$, we define the projection $p:P \to M$ by $p:=\pi \circ \pi_2$.
Under this identification, the exchange map $s: Y \times_M Y_k \to Y_k \times_M
Y$ becomes an involution of $P$
which lifts $k$,
\begin{equation}
\alxy{P \ar[d]_{p} \ar[r]^{s} & P \ar[d]^{p} \\ M \ar[r]_{k} & M\text{.}}
\end{equation}  
 The dual 1-isomorphism
$\mathcal{A}^{*}$ has by definition the line bundle $s^{*}A$ over $P$. Now, similarly as for the pullback of $\pi:Y \to M$ we denote the pullback
of $p:P \to M$ by $p_k: P_k \to M$ and choose $P_k:=P$ and $p_k:=k
\circ p$. This way, the pullback 1-isomorphism $k^{*}\mathcal{A}$ has
the line bundle $A$ over $P$. Again by Theorem \ref{th1}, we assume that
the 2-isomorphism $\varphi$ can be represented by a triple $(P,\id_P,\varphi_P)$ with an isomorphism $\varphi_P: A \to
s^{*}A$ of line bundles over $P$ satisfying the compatibility axiom \hyperlink{2M}{(2M)} with the
isomorphism $\alpha$ of $\mathcal{A}$:
\begin{equation}
\label{28}
\alxydim{@C=1.4cm@R=1.4cm}{L \otimes \zeta_2^{*}A \ar[d]_{\id \otimes \zeta_2^{*}\varphi_P} \ar[r]^{\alpha} & \zeta_1^{*}A \otimes L \ar[d]^{\zeta_1^{*}\varphi_P \otimes \id} \\ L \otimes \zeta_2^{*}s^{*}A \ar[r]_{s^{*}\alpha}
& \zeta_1^{*}s^{*}A \otimes L}
\end{equation}
The dual 2-isomorphism $\varphi^{*}$ is given by $(P,\id_P,s^{*}\varphi_P)$, and
the equation $\varphi = k^{*}\varphi^{*-1}$
becomes $\varphi_P = s^{*}\varphi^{-1}_P$. So,
$\varphi_P$ is an $s$-equivariant structure on $A$. This is exactly the original
definition \cite{schreiber1}: a stable isomorphism $\mathcal{A}: k^{*}\mathcal{G} \to \mathcal{G}^{*}$, whose line bundle $A$ is equipped with an $s$-equivariant structure which
is  compatible with the isomorphism $\alpha$ of $\mathcal{A}$ in the sense
of the commutativity of diagram (\ref{28}).

\medskip

Defining a Jandl structure in terms of 1- and 2-morphisms has -- just like
for gerbe modules -- several advantages. For example, it is easy to see
that Jandl structures are compatible with pullbacks along equivariant maps, tensor products and
duals of bundle gerbes. Furthermore, we have an obvious definition of morphisms
between Jandl structures, which induces exactly the notion of equivalent
Jandl structures we introduced in \cite{schreiber1}.
\begin{definition}
A morphism $\beta:\mathcal{J} \to \mathcal{J}'$ between Jandl structures $\mathcal{J}=(k,\mathcal{A},\varphi)$ and $\mathcal{J}'=(k,\mathcal{A}',\varphi')$
on the same bundle gerbe $\mathcal{G}$ over $M$
with the same involution $k$ is a 2-morphism 
\begin{equation*}
\beta\lmaps \mathcal{A} \Longrightarrow
\mathcal{A}'
\end{equation*}
 which commutes with $\varphi$
and $\varphi'$ in the sense that the diagram
\begin{equation*}
\alxydim{@C=1.5cm@R=1cm}{\mathcal{A} \ar@{=>}[r]^{\varphi}
\ar@{=>}[d]_{\beta} & k^{*}\mathcal{A}^{*} \ar@{=>}[d]^{k^{*}\beta^{*}}
\\ \mathcal{A}' \ar@{=>}[r]_{\varphi'}
& k^{*}\mathcal{A}'^{*}}
\end{equation*}
of 2-morphisms is commutative.
\end{definition}
Since $\mathcal{A}$ is invertible, every morphism of Jandl structures
is invertible. We may thus speak of a groupoid $\mathfrak{Jdl}(\mathcal{G},k)$ of Jandl structures on the bundle gerbe $\mathcal{G}$ with involution $k$. The skeleton
of this groupoid has been classified \cite{schreiber1}: it forms a torsor over the group of flat $k$-equivariant
line bundles over $M$. The following proposition relates these groupoids
 of Jandl structures
on isomorphic bundle gerbes on the same space with the same involution. This relation is a new result, coming and benefiting
very much from the 2-categorial structure of bundle gerbes we have developed.

\begin{proposition}
\label{prop3}
Any 1-isomorphism $\mathcal{B}:\mathcal{G} \to \mathcal{G}'$ induces an equivalence
of groupoids
\begin{equation*}
J_{\mathcal{B}}\lmaps \mathfrak{Jdl}(\mathcal{G}',k) \lto \mathfrak{Jdl}(\mathcal{G},k)
\end{equation*}
with the following properties:
\begin{itemize}
\item[a)]
any 2-morphism $\beta: \mathcal{B}
\Rightarrow \mathcal{B}'$ induces a natural equivalence $J_{\mathcal{B}} \cong
J_{\mathcal{B}}'$. 
\item[b)]
there is a  natural equivalence $J_{\id_{\mathcal{G}}} \cong \id_{\mathfrak{Jdl}(\mathcal{G},k)}$.
\item[c)]
it respects the composition of 1-morphisms in the sense that
\begin{equation*}
J_{\mathcal{B}' \circ \mathcal{B}} = J_{\mathcal{B}} \circ J_{\mathcal{B}'}
\text{.}
\end{equation*}
\end{itemize}

\end{proposition}

\proof
The functor $J_{\mathcal{B}}$ sends a Jandl structure $(k,\mathcal{A},\varphi)$ on $\mathcal{G}'$
to the triple $(k,\mathcal{A}',\varphi')$ with the same involution $k$, the 1-isomorphism
\begin{equation}
\label{105}
\mathcal{A}' := \mathcal{B}^{*} \circ  \mathcal{A} \circ k^{*}\mathcal{B}
\lmaps k^{*}\mathcal{G} \lto \mathcal{G}^{*}
\end{equation}
and the 2-isomorphism
\begin{equation}
\label{79}
\alxydim{@C=0.6cm}{k^{*}\mathcal{A}' \ar@{=}[r] & k^{*}\mathcal{B}^{*}\circ
k^{*}\mathcal{A}
\circ
\mathcal{B} \ar@{=>}^{\id_{k^{*}\mathcal{B}^{*}} \circ \varphi \circ \id_{\mathcal{B}}}[d] & \\ & k^{*}\mathcal{B}^{*}\circ
\mathcal{A}^{*}
\circ
\mathcal{B}  \ar@{=}[r] & k^{*}\mathcal{A}'^{*}}
\end{equation}
where we use equation (\ref{82}).
The following calculation shows that $(k,\mathcal{A}',\varphi')$ is a Jandl
structure:
\begin{eqnarray}
k^{*}\varphi'^{*} &\stackrel{\text{def}}{=}&  k^{*}(\id_{k^{*}\mathcal{B}^{*}} \circ \varphi \circ \id_{\mathcal{B}})^{*}
\nonumber \\
 &\stackrel{\text{(\ref{82})}}{=}& \id_{k^{*}\mathcal{B}^{*}}\circ k^{*}\varphi^{*} \circ \id_{\mathcal{B}}\nonumber \\
&\stackrel{}{=}& \id_{\mathcal{B}}\circ \varphi^{-1} \circ \id_{\mathcal{B}^{*}}
\nonumber \\ &\stackrel{\text{def}}{=}& \varphi'^{-1}\text{.}
\end{eqnarray}
A morphism $\beta$ of Jandl structures on $\mathcal{G}'$
is sent to the morphism
\begin{equation}
J_{\mathcal{B}}(\beta) := \id_{\mathcal{B}^{*}} \circ \beta \circ \id_{k^{*}\mathcal{B}}
\end{equation}
of the respective Jandl structures on $\mathcal{G}'$.
The two axioms of the composition functor $\circ$ from Lemma \ref{lem5} show
that the composition of morphisms of Jandl structures is respected, so that
$J_{\mathcal{B}}$ is a functor.
It is an equivalence because $J_{\mathcal{B}^{-1}}$ is an inverse functor, where the natural equivalences $J_{\mathcal{B}^{-1}}
\circ J_{\mathcal{B}} \cong \id$ and $J_{\mathcal{B}}
\circ J_{\mathcal{B}^{-1}} \cong \id$ use the 2-isomorphisms $i_r$ and $i_l$
from section \ref{sec6} associated to the inverse 1-morphism $\mathcal{B}^{-1}$.

To prove a), let $\beta: \mathcal{B} \Rightarrow \mathcal{B}'$ be a 2-morphism.
We define
the natural equivalence $J_{\mathcal{B}} \cong J_{\mathcal{B}'}$, which
is a collection
of morphisms $\beta_{\mathcal{J}}: J_{\mathcal{B}}(\mathcal{J}) \to
J_{\mathcal{B}'}(\mathcal{J})$ of Jandl structures on $\mathcal{G}$ for any
Jandl structure $\mathcal{J}$ on $\mathcal{G}'$ by
\begin{equation}
\beta_{\mathcal{J}} :=
\beta^{*}\circ \id_{\mathcal{A}} \circ k^{*}\beta\text{.}
\end{equation}
 This defines indeed a
morphism of Jandl structures and makes the naturality square
\begin{equation}
\alxydim{@C=1.7cm@R=1.4cm}{J_{\mathcal{B}}(\mathcal{J}) \ar[r]^{\beta_{\mathcal{J}}} \ar[d]_{J_{\mathcal{B}}(\beta)}
& J_{\mathcal{B}'} (\mathcal{J}) \ar[d]^{J_{\mathcal{B}'}(\beta)}\\ J_{\mathcal{B}}(\mathcal{J}')
\ar[r]_{\beta_{\mathcal{J}'}} & J_{\mathcal{B}'}(\mathcal{J}')}
\end{equation}
commutative. The natural equivalence for b) uses the 2-isomorphisms $\lambda_{\mathcal{A}}$
and $\rho_{\mathcal{A}}$ of the 2-category $\mathfrak{BGrb}(M)$ and the fact
that $\id_{\mathcal{G}}^{*}=\id_{\mathcal{G}^{*}}$. Finally,  c) follows  from the definition of $J_{\mathcal{B}}$
and the fact that the duality functor $()^{*}$ respects the composition of
1-morphisms, see (\ref{82}).  
\endofproof

It is worthwhile to consider a Jandl structure $\mathcal{J}=(k,\mathcal{A},\varphi)$ over a trivial bundle gerbe
$\mathcal{I}_{\rho}$. By definition, this is a 1-isomorphism
\begin{equation}
\mathcal{A} \lmaps \mathcal{I}_{k^{*}\rho} \lto \mathcal{I}_{-\rho}
\end{equation} 
and a 2-isomorphism $\varphi :k^{*}\mathcal{A} \Rightarrow \mathcal{A}^{*}$
satisfying $\varphi = k^{*}\varphi^{*-1}$. Now we apply the functor $\mathrm{Bun}$ and obtain
a line bundle $\hat R := \mathrm{Bun}(\mathcal{A})$ over $M$ of curvature
$-(\rho +k^{*} \rho)$ and an isomorphism
$\hat\varphi:= \mathrm{Bun}(\varphi): k^{*}\hat R \to \hat R$
 of line bundles over $M$ which satisfies $\hat \varphi = k^{*}\hat\varphi^{-1}$,
 summarizing: a $k$-equivariant line bundle. So, the functor $\mathrm{Bun}$ induces an equivalence of groupoids
 \begin{equation}
\mathrm{Bun}^{k}_{\rho} \lmaps \mathfrak{Jdl}(\mathcal{\mathcal{I}_{\rho}},k) \lto \mathfrak{LBun}_{-(\rho+k^{*}\rho)}^{k}(M)
\end{equation}
between the groupoid of  Jandl structures  on $\mathcal{I}_{\rho}$ with involution
$k$ and the groupoid of $k$-equivariant
line bundles over $M$ with curvature $-(\rho+k^{*}\rho)$. In particular,
if $\mathcal{G}$ is a bundle gerbe over $M$ and $\mathcal{T}:\mathcal{G}
\to \mathcal{I}_{\rho}$  a trivialization, we obtain a functor
\begin{equation}
\label{80}
\alxydim{@C=1.7cm}{\mathfrak{Jdl}(\mathcal{G},k) \ar[r]^-{J_{\mathcal{T}^{-1}}}  &  \mathfrak{Jdl}(\mathcal{I}_{\rho},k)
\ar[r]^-{\mathrm{Bun}^k_{\rho}} & \mathfrak{LBun}_{-(\rho+k^{*}\rho)}^{k}(M)}
\end{equation}
converting a Jandl structure on the bundle gerbe $\mathcal{G}$ into a $k$-equivariant
line bundle over $M$. It becomes obvious that the existence of a Jandl
structure with involution $k$ on the trivial bundle gerbe $\mathcal{I}_{\rho}$ constraints the
2-form $\rho$: as the curvature of a line bundle, the 2-form $-(\rho+k^{*}\rho)$ has to be closed and to have integer periods. 

\medskip

Let us now explain how Jandl structures enter in the definition of holonomy
around unoriented surfaces, and how we can take further advantage of the 2-categorial
formalism. We have learned before that to incorporate surfaces with boundary
we had to do two steps: we first specified additional structure, a D-brane
of the bundle gerbe $\mathcal{G}$,
and then specified which maps $\phi: \Sigma \to M$ are compatible with this
additional structure: those who send the boundary of $\Sigma$
into the support of the D-brane. To discuss unoriented surfaces (without
boundary), we also
do these two steps: the additional structure we choose here is a
Jandl structure $\mathcal{J}=(k,\mathcal{A},\varphi)$ on the bundle gerbe
$\mathcal{G}$. To describe the space of maps we want to consider, we have to introduce
the following geometric structures \cite{schreiber1}:

\begin{itemize}
\item 
For any (unoriented) closed surface $\Sigma$ there
is an oriented two-fold covering $\mathrm{pr}:\hat\Sigma \to
\Sigma$. It is unique up to orientation-preserving
diffeomorphisms and it is connected if and
only if $\Sigma$ is not orientable. It has a canonical,
orientation-reversing involution
$\sigma$, which permutes the sheets
and preserves the fibres. We call this two-fold covering the
orientation  covering of $\Sigma$. 

\item A fundamental domain  of
$\Sigma$ in $\hat\Sigma$ is a submanifold $F$ of $\hat\Sigma$ with
( possibly only piecewise smooth) boundary, such that 
\begin{equation}
F \cup \sigma(F) = \hat\Sigma
\quad\text{ and }\quad
F \cap \sigma(F) = \partial F\text{.}
\end{equation}
A key
observation is that the involution
$\sigma$ restricts to an orientation-{\it
preserving}  involution on $\partial F\subset
 \hat\Sigma$. Accordingly, the quotient $\overline{\partial
 F}$ is an oriented closed
 1-dimensional submanifold of $\Sigma$. 

\end{itemize}
Now, given a closed surface $\Sigma$, we consider maps $\hat\phi:\hat\Sigma
\to M$  from the orientation covering
$\hat\Sigma$ to $M$, which are equivariant with respect to the
two involutions on $\hat\Sigma$ and $M$,
i.e. the diagram 
\begin{equation*}
\alxydim{@C1.3cm@R1.3cm}{\hat\Sigma \ar[r]^{\hat \phi}
\ar[d]_{\sigma} & M \ar[d]^{k}
\\ \hat\Sigma \ar[r]_{\hat\phi}
& M}
\end{equation*} 
has to be commutative.

\begin{definition}
\label{def3}
Let $\mathcal{J}$ be a Jandl structure on a bundle gerbe $\mathcal{G}$
over $M$, and let $\hat\phi: \hat\Sigma \to M$ be an equivariant
smooth map. For a trivialization
\begin{equation*}
\mathcal{T}\lmaps \hat\phi^{*}\mathcal{G} \lto \mathcal{I}_{\rho}
\end{equation*}
let $\hat R$ be the  $\sigma$-equivariant line bundle over $\hat\Sigma$ determined
by the functor 
\begin{equation}
\mathrm{Bun}^{\rho}_\sigma \circ J_{\mathcal{T}^{-1}}: \mathfrak{Jdl}(\hat\phi^{*}\mathcal{G},\sigma)
\lto \mathfrak{LBun}_{-(\rho+\sigma^{*}\rho)}^{\sigma}(\hat \Sigma) 
\end{equation}
from (\ref{80}). 
In turn, $\hat R$ defines a line bundle $R$ over $\Sigma$. Choose any fundamental
domain $F$ of $\Sigma$. Then, the holonomy of $\mathcal{G}$ with Jandl structure
$\mathcal{J}$ around $\hat\phi$ is defined as
\begin{equation*}
\mathrm{hol}_{\mathcal{G},\mathcal{J}}(\hat\phi) := \mathrm{exp} \left (
\mathrm{i} \int_{F}\rho 
\right ) \cdot \mathrm{hol}_R(\overline{\partial F})\text{.}
\end{equation*}
\end{definition}

Definition \ref{def3} is a generalization of Definition \ref{def6} of holonomy
around an oriented surface: for an orientable surface $\Sigma$ and 
\textit{any} choice of an orientation, they coincide \cite{schreiber1}. To show that
Definition \ref{def3} does not depend on the choice of the trivialization
$\mathcal{T}$, we combine all the collected tools. Let  $\mathcal{T}': \hat\phi^{*}\mathcal{G}
\to \mathcal{I}_{\rho'}$ be any other trivialization. We consider the 1-isomorphism \begin{equation}
\mathcal{B}\;:=\; \mathcal{T} \circ \mathcal{T}'^{-1} \lmaps \mathcal{I}_{\rho'} \lto \mathcal{I}_{\rho}
\end{equation}
and
the corresponding line bundle $T:=\mathrm{Bun}(\mathcal{B})$. To compare
the two $\sigma$-equivariant line bundles $\hat R$ and $\hat R'$ corresponding
to the two trivializations, we first compare the Jandl structures $J_{\mathcal{T}^{-1}}(\mathcal{J})$
on $\mathcal{I}_{\rho}$ and $J_{\mathcal{T}'^{-1}}(\mathcal{J})$ on $\mathcal{I}_{\rho'}$.
By Proposition \ref{prop3} a), b) and c), there exists an isomorphism 
\begin{equation}
J_{\mathcal{T}'^{-1}}(\mathcal{J}) \cong J_{\mathcal{B}} (J_{\mathcal{T}^{-1}}(\mathcal{J}))
\end{equation}  
of Jandl structures on $\mathcal{I}_{\rho}$. 
By definition of the functor $J_{\mathcal{B}}$, this isomorphism is a 2-isomorphism
\begin{equation}
\mathcal{A}' \cong\mathcal{B}^{*} \circ \mathcal{A} \circ
\sigma^{*}\mathcal{B}\text{,}
\end{equation}
where $\mathcal{A}$ is the 1-morphism of $J_{\mathcal{T}^{-1}}(\mathcal{J})$
and $\mathcal{A}'$ is the 1-morphism of $J_{\mathcal{T}'^{-1}}(\mathcal{J})$. Now we apply the functor $\mathrm{Bun}$
and
obtain an isomorphism
\begin{equation}
\hat R' \cong  T \otimes \hat R \otimes \sigma^{*}T 
\end{equation}
of $\sigma$-equivariant line bundles over $\hat \Sigma$, where $\hat Q:= \sigma^{*}T \otimes T$ has the canonical $\sigma$-equivariant
structure by exchanging the tensor factors. Thus, we have isomorphic line bundles
\begin{equation}
R'\cong R \otimes Q
\end{equation}  
over $\Sigma$. Notice that the holonomy of the line bundle
$Q$ is
\begin{equation}
\mathrm{hol}_Q(\overline{\partial F}) = \mathrm{hol}_{T}(\partial F) = \exp
\left ( \mathrm{i} \int_F \rho - \rho' \right ) 
\end{equation} 
This shows
\begin{eqnarray}
\mathrm{exp} \left (
\mathrm{i} \int_{F}\rho' 
\right ) \cdot \mathrm{hol}_{R'}(\overline{\partial F}) &=& \mathrm{exp} \left (
\mathrm{i} \int_{F}\rho'
\right ) \cdot \mathrm{hol}_{Q}(\overline{\partial F}) \cdot \mathrm{hol}_{R }(\overline{\partial F}) 
\nonumber\\ &=& \mathrm{exp} \left (
\mathrm{i} \int_{F}\rho
\right ) \cdot \mathrm{hol}_{R}(\overline{\partial F}) 
\end{eqnarray}
so that Definition \ref{def3} does not depend on the choice of the trivialization.
In \cite{schreiber1} we have deduced from the equation $\mathrm{curv}(\hat R) = -(\rho
+ \sigma^{*}\rho)$ that it is also independent of the choice of the fundamental
domain. 

\medskip

Unoriented surface holonomy, defined in terms of Jandl structures on bundle
gerbes, provides a candidate for the Wess-Zumino term in two-dimensional
conformal field theory for unoriented worldsheets,  as they appear in type
I string theories. Following the examples of $M=SU(2)$ and $M=SO(3)$
we give in \cite{schreiber1},
we reproduce results known from other approaches. This indicates, that a
bundle gerbe with Jandl structure, together with a metric, is the background field for unoriented WZW
models. In this setup, Proposition \ref{prop3} assures, that -- just like
for oriented WZW models -- only the isomorphism class
of the bundle gerbe is relevant.

\newcommand{\etalchar}[1]{$^{#1}$}

\end{document}